\documentclass[11pt]{article}
\usepackage{latexsym,amsfonts,amssymb,amsmath,amsthm}
\usepackage{graphicx}

\usepackage[usenames,dvipsnames]{color}
\usepackage{ulem}

\begin{document}
\setlength{\baselineskip}{16pt}

\parindent 0.5cm
\evensidemargin 0cm \oddsidemargin 0cm \topmargin 0cm \textheight
22cm \textwidth 16cm \footskip 2cm \headsep 0cm

\newtheorem{theorem}{Theorem}[section]
\newtheorem{lemma}{Lemma}[section]
\newtheorem{proposition}{Proposition}[section]
\newtheorem{definition}{Definition}[section]
\newtheorem{example}{Example}[section]
\newtheorem{corollary}{Corollary}[section]

\newtheorem{remark}{Remark}[section]

\numberwithin{equation}{section}

\def\p{\partial}
\def\I{\textit}
\def\R{\mathbb R}
\def\C{\mathbb C}
\def\u{\underline}
\def\l{\lambda}
\def\a{\alpha}
\def\O{\Omega}
\def\e{\epsilon}
\def\ls{\lambda^*}
\def\D{\displaystyle}
\def\wyx{ \frac{w(y,t)}{w(x,t)}}
\def\imp{\Rightarrow}
\def\tE{\tilde E}
\def\tX{\tilde X}
\def\tH{\tilde H}
\def\tu{\tilde u}
\def\d{\mathcal D}
\def\aa{\mathcal A}
\def\DH{\mathcal D(\tH)}
\def\bE{\bar E}
\def\bH{\bar H}
\def\M{\mathcal M}
\renewcommand{\labelenumi}{(\arabic{enumi})}

\def\disp{\displaystyle}
\def\undertex#1{$\underline{\hbox{#1}}$}
\def\card{\mathop{\hbox{card}}}
\def\sgn{\mathop{\hbox{sgn}}}
\def\exp{\mathop{\hbox{exp}}}
\def\OFP{(\Omega,{\cal F},\PP)}
\newcommand\JM{Mierczy\'nski}
\newcommand\RR{\ensuremath{\mathbb{R}}}
\newcommand\CC{\ensuremath{\mathbb{C}}}
\newcommand\QQ{\ensuremath{\mathbb{Q}}}
\newcommand\ZZ{\ensuremath{\mathbb{Z}}}
\newcommand\NN{\ensuremath{\mathbb{N}}}
\newcommand\PP{\ensuremath{\mathbb{P}}}
\newcommand\abs[1]{\ensuremath{\lvert#1\rvert}}

\newcommand\normf[1]{\ensuremath{\lVert#1\rVert_{f}}}
\newcommand\normfRb[1]{\ensuremath{\lVert#1\rVert_{f,R_b}}}
\newcommand\normfRbone[1]{\ensuremath{\lVert#1\rVert_{f, R_{b_1}}}}
\newcommand\normfRbtwo[1]{\ensuremath{\lVert#1\rVert_{f,R_{b_2}}}}
\newcommand\normtwo[1]{\ensuremath{\lVert#1\rVert_{2}}}
\newcommand\norminfty[1]{\ensuremath{\lVert#1\rVert_{\infty}}}

\title{Liouville Type Property and Spreading Speeds of KPP Equations
in Periodic Media with Localized Spatial  Inhomogeneity}

\author{
Liang Kong\\
Department of Mathematical Sciences\\
   University of Illinois Springfield\\
   Springfield, IL 62703, USA\\
   \\
 and\\
 \\
  Wenxian Shen\thanks{Partially supported by NSF grant DMS--0907752}\\
Department of Mathematics and Statistics\\
Auburn University\\
Auburn, AL 36849, U.S.A. }
\date{}
\maketitle

\noindent{\bf Abstract.}
The current paper is devoted to the study of semilinear dispersal evolution equations of the form
$$
u_t(t,x)=(\mathcal{A}u)(t,x)+u(t,x)f(t,x,u(t,x)),\quad x\in\mathcal{H},
$$
where $\mathcal{H}=\RR^N$ or $\ZZ^N$, $\mathcal{A}$ is a random dispersal operator or nonlocal dispersal
operator in the case $\mathcal{H}=\RR^N$ and is a discrete dispersal operator in the case $\mathcal{H}=\ZZ^N$,
and $f$ is periodic in $t$, asymptotically periodic in $x$ (i.e. $f(t,x,u)-f_0(t,x,u)$ converges to $0$ as $\|x\|\to\infty$ for some
time and space periodic function $f_0(t,x,u)$), and is of KPP type in $u$. It is proved that Liouville type
property for such equations holds, that is, time periodic  strictly positive solutions are unique. It is also proved that
if $u\equiv 0$
is a linearly unstable solution to the time and space periodic limit equation of such an equation, then it
 has a unique stable time periodic strictly  positive
solution and has a spatial spreading speed in every direction.

\medskip


\medskip


\section{Introduction}

In this paper, we consider the existence, uniqueness, and stability of time periodic positive solutions and spatial spreading speeds of KPP type
dispersal evolution equations in periodic media with  localized spatial inhomogeneity.
Our model equations are of the form,
\begin{equation}
\label{periodic-eq}
u_t(t,x)=\mathcal{A} u+uf_0(t,x,u),\quad x\in\mathcal{H},
\end{equation}
where $\mathcal{H}=\RR^N$ or $\ZZ^N$;  in the case $\mathcal{H}=\RR^N$,
$\mathcal{A}u=\Delta u$ or $(\mathcal{A}u)(t,x)=\int_{\RR^N}\kappa(y-x)u(t,y)dy-u(t,x)$
($\kappa(\cdot)$ is a smooth non-negative
 convolution kernel supported on a ball centered at the origin and $\int_{\RR^N}\kappa(z)dz=1$),
and in the case $\mathcal{H}=\ZZ^N$, $(\mathcal{A}u)(t,j)= \sum_{k\in K}a_k(u(t,j+k)-u(t,j))$ ($a_k>0$ and $K=\{k\in\mathcal{H}\,|\,\|k\|=1\}$); and $f_0(t+T,x,u)=f_0(t,x+p_i{\bf e_i},u)=f(t,x,u)$
($T\in\RR$ and $p_i\in\mathcal{H}$  are given constants),
and $\p_u f(t,x,u)<0$ for $u\ge 0$, $f(t,x,u)<0$ for $u\gg 1$.

Among others, equation \eqref{periodic-eq} is used to model the evolution of population density of a species. The case that $\mathcal{H}=\RR^N$ and $\mathcal{A}u=\Delta u$ indicates
that the environment of the underlying model problem is not patchy and the internal interaction  of the organisms is
 random and local
(i.e. the organisms move randomly between the adjacent spatial locations, such $\mathcal{A}$ is referred to as  a {\it random dispersal operator}) (see
\cite{ArWe1}, \cite{ArWe2}, \cite{CaCo},  \cite{Fif1}, \cite{FiPe}, \cite{Fisher}, \cite{KPP},
\cite{Mur},  \cite{ShKa},
\cite{Ske}, \cite{Wei1}, \cite{Wei2}, \cite{Zha2}, etc., for the application in this case).
 If the environment of the underlying model problem is not patchy
and the internal interaction of the organisms is nonlocal, $(\mathcal{A}u)(t,x)=\int_{\RR^N}\kappa(y-x)u(t,y)dy-u(t,x)$ is often adopted (such $\mathcal{A}$ is referred to  as a {\it nonlocal dispersal
operator})
(see \cite{BaZh}, \cite{CaCoLoRy},  \cite{ChChRo}, \cite{CoCoElMa}, \cite{Fif2}, \cite{GrHiHuMiVi}, \cite{HuMaMiVi}, \cite{KaLoSh},
 etc.).
The case that $\mathcal{H}=\ZZ^N$ and $(\mathcal{A}u)(t,j)= \sum_{k\in K}a_k(u(t,j+k)-u(t,j))$ (which is referred to as a {\it discrete dispersal operator})
arises when modeling the population dynamics of species living in  patchy environments
(see \cite{Fif1},  \cite{MaWuZo}, \cite{Mur},  \cite{ShKa}, \cite{ShSw}, \cite{Wei1}, \cite{Wei2}, \cite{WuZo},  etc.).
The periodicity of $f_0(t,x,u)$ in $t$ and $x$ reflects the periodicity of the environment. In literature,
equation \eqref{periodic-eq} is called  Fisher or  KPP type due to the pioneering works of Fisher \cite{Fisher} and Kolmogorov, Petrowsky, Piscunov
\cite{KPP}
 on the following special case of
\eqref{periodic-eq},
\begin{equation}
\label{fisher-kpp-eq}
u_t=u_{xx}+u(1-u),\quad x\in\RR.
\end{equation}

Central problems about \eqref{periodic-eq} include the existence, uniqueness, and stability of time and space periodic positive solutions and spatial spreading speeds.
Such problems   have been extensively studied (see \cite{ArWe1}-\cite{BeHaRo}, \cite{ChFuGu}-\cite{ChGu2}, \cite{CoDaMa}, \cite{CoDu}, \cite{FrGa}, \cite{GuHa}-\cite{GuWu2}, \cite{HeShZh},
\cite{HuZi1}, \cite{Kam}, \cite{LiSuWa}-\cite{LiYiZh},  \cite{Nad0}-\cite{NoXi}, \cite{RaSh}-\cite{ShZh3}, \cite{Uch}-\cite{Wei2}, etc.).
It is known that time and space periodic positive solutions of \eqref{periodic-eq} (if exist)  are unique, which is referred to as  the  {\it Liouville type property} for \eqref{periodic-eq}.
If $u\equiv 0$ is linearly unstable with respect to spatially periodic perturbations, then \eqref{periodic-eq} has a unique stable time and space periodic positive solution
$u_0^*(t,x)$  and for any $\xi\in \RR^N$ with $\|\xi\|=1$,
\eqref{periodic-eq} has a spreading speed $c_0^*(\xi)$ in the direction of $\xi$ (see section 2.4 for detail).

The aim of the current paper is to deal with the extensions of the above results for \eqref{periodic-eq} to KPP type equations in periodic media with spatially localized
inhomogeneity. We  consider
\begin{equation}
\label{main-eq}
u_t=\mathcal{A}u+uf(t,x,u),\quad x\in\mathcal{H},
\end{equation}
where $\mathcal{A}$ and $\mathcal{H}$ are as in \eqref{periodic-eq}, $\p_u f(t,x,u)<0$ for $u\ge 0$, $f(t,x,u)<0$ for $u\gg 1$,
$f(t+T,x,u)=f(t,x,u)$, and  $|f(t,x,u)-f_0(t,x,u)|\to 0$ as $\|x\|\to\infty$ uniformly in $(t,u)$ on bounded sets
($f_0(t,x,u)$ is as in \eqref{periodic-eq}). We show that localized inhomogeneity does not destroy the existence and uniqueness
of time periodic positive solutions and it neither slows down nor speeds up the spatial spreading speeds. More precisely, we prove

\medskip

\noindent $\bullet$ (Liouville type property or uniqueness of time periodic strictly positive solutions)
{\it Time periodic strictly positive solutions of \eqref{main-eq} (if exist) are unique (see Theorem \ref{positive-solution-thm}(1)).}
\medskip

\noindent $\bullet$ (Stability of time periodic strictly  positive solutions)
{\it  If \eqref{main-eq} has a time periodic strictly positive solution
$u^*(t,x)$, then it is asymptotically stable (see Theorem \ref{positive-solution-thm}(2)).}

\medskip

\noindent $\bullet$ (Existence of time periodic strictly positive solutions) {\it If $u=0$ is a linearly unstable solution of \eqref{periodic-eq} with respect to
periodic perturbations, then \eqref{main-eq} has a time periodic strictly positive solution $u^*(t,x)$ (see Theorem \ref{positive-solution-thm}(3)).}

\medskip

\noindent $\bullet$ (Tail property of time periodic strictly positive solutions) {\it If $u=0$ is a linearly unstable solution of \eqref{periodic-eq} with respect to
periodic perturbations, then  $u^*(t,x)-u_0^*(t,x)\to 0$ as $\|x\|\to \infty$ uniformly in $t$ (see Theorem \ref{positive-solution-thm}(4)).
}
\medskip

\noindent $\bullet$ (Spatial spreading speeds) {\it If $u=0$ is a linearly unstable solution of \eqref{periodic-eq} with respect to
periodic perturbations, then for each $\xi\in\RR^N$ with $\|\xi\|=1$, $c_0^*(\xi)$ is the spreading speed of \eqref{main-eq}
in the direction of $\xi$ (see Theorem \ref{spreading-speed-thm}).

\medskip

\noindent $\bullet$ (Spreading features of spreading speeds) The spreading speeds of \eqref{main-eq} are of important spreading features
(see Theorem \ref{spreading-feature-thm} for detail).
}

\medskip

It should be pointed out that the Liouville type property for \eqref{main-eq} in the case that
$\mathcal{H}=\RR^N$, $\mathcal{A}=\Delta u$, and $f(t,x,u)=f(x,u)$ has been proved in \cite{BeHaRo}.
Existence, uniqueness, and stability of time independent positive solutions and spreading speeds
for \eqref{main-eq} in the case $f(t,x,u)=f(x,u)$ and $f_0(t,x,u)=f_0(u)$ have been proved in \cite{KoSh}.
However, many techniques developed in \cite{BeHaRo} and \cite{KoSh} are difficult to apply to \eqref{main-eq}. Several important
new techniques are developed in the current paper. Both the results and techniques established in
the current paper can be extended to more general cases (say,  cases that $\mathcal{A}$ is some linear combination of
random and nonlocal dispersal operators and/or  $f(t,x,u)$ is almost periodic
in $t$ and asymptotically  periodic in $x$).

It should also be pointed out that, if $u=0$ is a linearly unstable solution of \eqref{periodic-eq} with respect to
periodic perturbations, then \eqref{periodic-eq} has traveling wave solutions connecting
$0$ and $u_0^*(\cdot,\cdot)$ and propagating in the direction of $\xi$ with speed $c> c_0^*(\xi)$ for any
$\xi\in S^{N-1}$. But \eqref{main-eq} may have no traveling wave solutions connecting $0$ and $u^*(\cdot,\cdot)$
(see \cite{NoRoRyZl})
(hence, localized spatial inhomogeneity may prevent the existence of traveling wave solutions).

The rest of the paper is organized as follows. In section 2, we introduce
the standing notions, hypotheses, and definitions, and state the main results of the paper.
In section 3, we present some preliminary materials to be
used in the proofs of the main results. We study the existence, uniqueness, and stability of time periodic positive
solutions of \eqref{main-eq} in section 4. In section 5, we explore the spreading speeds of \eqref{main-eq}.

\section{Notions, Hypotheses, Definitions, and Main Results}

In this section, we first  introduce some standing notations,  hypotheses, and definitions.
We then state the main results of the paper.

\subsection{Notions, hypotheses and definitions}

In this subsection, we introduce standing notions, hypotheses, and definitions.
Throughout this subsection,  $\mathcal{H}=\RR^N$ or $\ZZ^N$ and $p_i\in\mathcal{H}$ with $p_i>0$
($i=1,2,\cdots,N$).

Let
\begin{equation}
\label{x-space-eq}
X=C_{\rm unif}^b(\mathcal{H},\RR):=\{u\in C(\mathcal{H},\RR)\,|\, u\,\, \text{is uniformly continuous and bounded on}\,\, \mathcal{H}\}
\end{equation}
with norm $\|u\|=\sup_{x\in\mathcal{H}}|u(x)|$,
\begin{equation}
\label{x-space-positive-eq}
X^+=\{u\in X\,|\, u(x)\ge 0\,\,\, \forall x\in\mathcal{H}\},
\end{equation}
and
\begin{equation}
\label{x-space-positive-positive-eq}
X^{++}=\{u\in X^+\,|\, \inf _{x\in\mathcal{H}}u(x)>0\}.
\end{equation}
Let
\begin{equation}
\label{x-periodic-space-eq}
X_p=\{u\in X\,|\, u(\cdot+p_i{\bf e_i})=u(\cdot)\},
\end{equation}
\begin{equation}
\label{x-periodic-space-positive-eq}
X_p^+=X^+\cap X_p,
\end{equation}
and
\begin{equation}
\label{x-periodic-space-positive-positive-eq}
X_p^{++}=X^{++}\cap X_p.
\end{equation}
 For given  $u,v\in X$, we define
\begin{equation}
\label{order-eq} u\leq v\,\, (u\geq v)\quad {\rm if}\,\, v-u\in
X^+\,\, (u-v\in X^+)
\end{equation}
and
\begin{equation}
\label{order-order-eq}
u\ll v\,\, (u\gg v)\quad {\rm if}\,\, v-u\in X^{++}\,\, (u-v\in X^{++}).
\end{equation}

Let $\mathcal{H}_i$ and $\mathcal{A}_i:\mathcal{D}(\mathcal{A}_i)\subset X\to X$ ($i=1,2,3)$ be defined by
\begin{equation}
\label{case1-eq}
\mathcal{H}_1=\RR^N,\,\, (\mathcal{A}_1 u)(x)=\Delta u(x) \quad \forall\,\, u\in \mathcal{D}(\mathcal{A}_1)(\subset X),
\end{equation}
\begin{equation}
\label{case2-eq}
\mathcal{H}_2=\RR^N,\,\, (\mathcal{A}_2u)(x)=\int_{\RR^N} \kappa(y-x)u(y)dy-u(x)\quad \forall\,\, u\in \mathcal{D}(\mathcal{A}_2)(= X),
\end{equation}
and
\begin{equation}
\label{case3-eq}
\mathcal{H}_3=\ZZ^N,\,\, (\mathcal{A}_3u)(j)=\sum_{k\in K}a_k(u(j+k)-u(j))\quad \forall u\in \mathcal{D}(\mathcal{A}_3)(=X).
\end{equation}

Let
\begin{equation}
\label{lift-x-space-eq}
\mathcal{X}_p=\{u\in C(\RR\times\mathcal{H},\RR)\,|\, u(t+T,x+p_i{\bf e_i})=u(t,x)\}
\end{equation}
with norm $\|u\|=\max_{t\in\RR,x\in\mathcal{H}}|u(t,x)|$.
For given $\xi\in S^{N-1}$ and $\mu\in\RR$, let $\mathcal{A}_{\xi,\mu}:\mathcal{D}(\mathcal{A}_{\xi,\mu})\subset\mathcal{X}_p\to\mathcal{X}_p$ be defined by
\begin{equation}
\label{operator-a-xi-mu-eq}
(\mathcal{A}_{\xi,\mu}u)(t,x)=\begin{cases}\Delta u(t,x)-2\mu \xi\cdot\nabla u(t,x)+\mu^2u(t,x)\,\, {\rm if}\,\, \mathcal{H}=\mathcal{H}_1,\,\,\mathcal{A}=\mathcal{A}_1\cr\cr
\int_{\RR^N} e^{-\mu (y-x)\cdot\xi}\kappa(y-x)u(t,y)dy-u(t,x)\,\,{\rm if}\,\, \mathcal{H}=\mathcal{H}_2,\,\,\mathcal{A}=\mathcal{A}_2\cr\cr
 \sum_{k\in K}a_k(e^{-\mu k\cdot\xi} u(t,j+k)-u(t,j))\,\,{\rm if}\,\, \mathcal{H}=\mathcal{H}_3,\,\,\mathcal{A}=\mathcal{A}_3
\end{cases}
\end{equation}
for $u\in\mathcal{D}(\mathcal{A}_{\xi,\mu})$.
Observe that
$$
\mathcal{A}_{\xi,0}=\mathcal{A}\quad \forall \,\,\xi\in S^{N-1}.
$$

 For any given $a\in \mathcal{X}_p$, $\xi\in S^{N-1}$, and $\mu\in\RR$,
let $\sigma(-\p_t+\mathcal{A}_{\xi,\mu}+a(\cdot,\cdot)\mathcal{I})$ be the spectrum of the operator
$-\p_t+\mathcal{A}_{\xi,\mu}+a(\cdot,\cdot)\mathcal{I}:\mathcal{D}(-\p_t+\mathcal{A}_{\xi,\mu}+a(\cdot,\cdot)\mathcal{I})\subset\mathcal{X}_p\to\mathcal{X}_p$,
$$
\big((-\p_t+\mathcal{A}_{\xi,\mu}+a(\cdot,\cdot)\mathcal{I})u\big)(t,x)
=-u_t(t,x)+(\mathcal{A}_{\xi,\mu}u)(t,x)+a(t,x)u(t,x).
$$
Let $\lambda_{\xi,\mu}(a)$ be defined by
\begin{equation}
\label{lambda-max-eq}
\lambda_{\xi,\mu}(a)=\sup\{{\rm Re}\lambda\,|\, \lambda\in\sigma (-\p_t+\mathcal{A}_{\xi,\mu}+a(\cdot,\cdot)\mathcal{I})\}.
\end{equation}
Observe that $\lambda_{\xi,0}(a)$ is independent of $\xi\in S^{N-1}$ and we may put
$$
\lambda(a)=\lambda_{\xi,0}(a).
$$

Consider \eqref{main-eq}. We introduce the following standing hypotheses.

\medskip
\noindent {\bf (H0)} {\it  $f(t+T,x,u)=f(t,x,u)$ for $(t,x,u)\in\RR\times\mathcal{H}\times\RR$ ($T>0$ is a given positive number);
$f(t,x,u)$ is $C^1$ in $t,u$ and $f(t,x,u)$, $f_t(t,x,u)$,  $f_u(t,x,u)$ are uniformly continuous in $(t,x,u)\in\RR\times\mathcal{H}\times E$ ($E$ is any bounded
 subset of $\RR$); $f(t,x,u)<0$ for all $t\in\RR$, $x\in\mathcal{H}$, and $u\ge M_0$ ($M_0>0$ is some given
constant); and $\inf_{t\in\RR, x\in\mathcal{H}}f_u(t,x,u)<0$ for all  $u\ge 0$.}
\medskip

\medskip
\noindent {\bf (H1)} {\it  $f_0(t,x,u)$ satisfies (H0), $f_0(t,x+p_i{\bf e_i},u)=f_0(t,x,u)$ for
$(t,x,u)\in\RR\times\mathcal{H}\times\RR$,  and $|f(t,x,u)-f_0(t,x,u)|\to 0$ as $\|x\|\to\infty$ uniformly in $(t,u)$ on bounded sets.}
\medskip

Throughout this section, we assume (H0).
By general semigroup theory (see \cite{Hen}, \cite{Paz}), for any $u_0\in X$, \eqref{main-eq} has a unique (local) solution
$u(t,\cdot;u_0)$ with $u(0,\cdot;u_0)=u_0(\cdot)$.
Furthermore, if $f(t,x+p_i{\bf e_i},u)=f(t,x,u)$ and   $u_0\in X_p$, then $u(t,\cdot;u_0)\in X_p$.
To indicate the dependence of $u(t,x;u_0)$ on $f$, we may write $u(t,x;u_0)$ as $u(t,x;u_0,f)$.

Let
\begin{equation}
\label{unit-sphere}
S^{N-1}=\{\xi\in\RR^N\,|\, \|\xi\|=1\}.
\end{equation}
For given $\xi\in S^{N-1}$ and $u\in X^+$,
we define
$$
\liminf_{x\cdot\xi\to -\infty}u(x)=\liminf_{r\to -\infty}\inf_{x\in\mathcal{H},x\cdot\xi\leq r}u(x).
$$
For given $u:[0,\infty)\times\mathcal{H}\to \RR$ and $c>0$, we define
$$
\liminf_{x\cdot\xi\leq ct, t\to\infty}u(t,x)=\liminf_{t\to\infty}\inf_{x\in\mathcal{H},x\cdot\xi\leq ct}u(t,x),
$$
$$
\limsup_{x\cdot\xi\geq ct,t\to\infty}u(t,x)=\limsup_{t\to\infty}\sup_{x\in\mathcal{H},x\cdot\xi\geq ct}u(t,x).
$$
The notions $\disp\limsup_{|x\cdot\xi|\leq ct,t\to\infty}u(t,x)$, $\disp\limsup_{|x\cdot\xi|\ge ct,t\to\infty}u(t,x)$,
$\disp\limsup_{\|x\|\leq ct,t\to\infty}u(t,x)$, and $\disp\limsup_{\|x\|\geq ct,t\to\infty}u(t,x)$ are defined similarly.
We define $X^+(\xi)$  by
\begin{equation}
\label{x-xi-space}
X^+(\xi)=\{u\in X^+\,|\, \liminf_{x\cdot\xi\to -\infty}u(x)>0,\quad u(x)=0\,\, {\rm for}\,\,  x\cdot\xi\gg 1\}.
\end{equation}

\begin{definition}[Time periodic strictly positive solution]
\label{positive-solution-def}
A solution $u(t,x)$ of \eqref{main-eq} on $t\in\RR$ is called a {\rm time periodic strictly positive solution}
if $u(t+T,x)=u(t,x)$ for $(t,x)\in\RR\times\mathcal{H}$ and $\inf_{(t,x)\in\RR\times\mathcal{H}}u(t,x)>0$.
\end{definition}

\begin{definition} [Spatial spreading speed]
\label{spreading-def}
For given $\xi\in S^{N-1}$, a real number $c^*(\xi)$
 is called the {\rm spatial spreading speed} of \eqref{main-eq} in the direction of $\xi$ if for any $u_0\in X^+(\xi)$,
$$
\liminf_{x\cdot\xi\leq ct, t\to \infty}u(t,x;u_0)>0\quad \forall c<c^*(\xi)
$$
and
$$
\limsup_{x\cdot\xi\geq ct, t\to\infty}u(t,x;u_0)=0\quad \forall c>c^*(\xi).
$$
\end{definition}

\subsection{Main Results}

In this subsection, we state the main results of this paper.

\begin{theorem} [Time periodic strictly  positive solutions]
\label{positive-solution-thm}
Consider \eqref{main-eq} and assume (H0).

\begin{itemize}

\item[(1)] $($Liouville type property or uniqueness$)$ If \eqref{main-eq} has a time periodic strictly positive solution, then it is unique.

\item[(2)] $($Stability$)$ Assume that $u^*(t,x)$ is a time periodic strictly positive solution of \eqref{main-eq}. Then
it is stable and for any $u_0\in X^{++}$,
$\lim_{t\to\infty}\|u(t,\cdot;u_0,f(\cdot+\tau,\cdot,\cdot))-u^*(t+\tau,
\cdot)\|_{X}=0$ uniformly in $\tau\in\RR$.

\item[(3)] $($Existence$)$
 Assume also (H1) and $\lambda(f_0(\cdot,\cdot,0))>0$.
 Then  \eqref{main-eq}
  has a unique time periodic strictly positive solution $u^*(t,x)$.

\item[(4)] $($Tail property$)$ Assume also (H1) and $\lambda(f_0(\cdot,\cdot,0))>0$. Then
$u^*(t,x)-u_0^*(t,x)\to 0$ as $\|x\|\to\infty$ uniformly in $t\in\RR$, where $u^*(t,x)$ is as in (3) and
$u_0^*(t,x)$ is the unique time and space periodic positive solution of \eqref{periodic-eq}
(see Proposition \ref{basic-periodic-solution} for the existence and uniqueness of $u_0^*(t,x)$).
\end{itemize}
\end{theorem}

\begin{theorem} [Existence of spreading speeds]
\label{spreading-speed-thm}
Consider \eqref{main-eq} and assume (H0) and (H1). If $\lambda(f_0(\cdot,\cdot,0))>0$, then
 for any given
$\xi\in S^{N-1}$, \eqref{main-eq}  has a spreading speed $c^*(\xi)$ in the direction of $\xi$. Moreover,
for any $u_0\in X^+(\xi)$,
\begin{equation}
\label{spreading-eq0}
\limsup_{x\cdot\xi\leq ct, t\to \infty}|u(t,x;u_0)-u^*(t, x)|=0\quad \forall c<c^*(\xi)
\end{equation}
and
 $$c^*(\xi)=c_0^*(\xi),
 $$
 where
$c_0^*(\xi)$
is the spatial spreading speeds of \eqref{periodic-eq} in the direction of $\xi$ (see Proposition \ref{basic-spreading-speed} for the existence
and characterization of $c_0^*(\xi)$).
\end{theorem}

\begin{theorem}[Spreading features of spreading speeds]
\label{spreading-feature-thm}
Consider \eqref{main-eq}. Assume (H0), (H1), and $\lambda(f_0(\cdot,\cdot,0))>0$.
 Then for any given $\xi\in S^{N-1}$, the following hold.
\begin{itemize}
\item[(1)] For each $ u_0\in X^+$ satisfying that  $u_0(x)=0$ for $x\in\mathcal{H}$ with $|x\cdot\xi|\gg 1$,
$$
\limsup_{|x\cdot\xi|\geq ct, t\to\infty}u(t,x;u_0)=0 \quad \forall\,\, c>\max\{c^*(\xi),c^*(-\xi)\}.
$$

\item[(2)]  For each $\sigma>0$,  $r>0$, and  $u_0\in X^+$ satisfying that $u_0(x)\geq\sigma$ for $x\in\mathcal{H}$ with $|x\cdot \xi|\leq r$,
$$
\limsup_{|x\cdot \xi|\leq ct,t\to\infty}|u(t,x;u_0)-u^*(t, x)|=0\quad \forall\,\, 0<c<\min\{c^*(\xi),c^*(-\xi)\}.
$$

\item[(3)] For each $u_0\in X^+$ satisfying that $u_0(x)=0$ for $x\in\mathcal{H}$ with $\|x\|\gg 1$,
$$
\limsup_{\|x\|\geq ct,t\to\infty}u(t,x;u_0)=0\quad \forall\,\, c>\sup_{\xi\in S^{N-1}}c^*(\xi).
$$

\item[(4)] For each $\sigma>0$, $r>0$, and  $u_0\in X^+$ satisfying that $u_0(x)\geq \sigma$ for $\|x\|\le r$,
$$
\limsup_{\|x\|\leq ct,t\to\infty}|u(t,x;u_0)-u^*(t, x)|=0\quad \forall\,\, 0<c<\inf_{\xi\in S^{N-1}}c^*(\xi).
$$
\end{itemize}
\end{theorem}

\section{Preliminary}

In this section, we present some preliminary materials to be used in later sections, including comparison principle for solutions of
\eqref{main-eq}; convergence of solutions of \eqref{main-eq} on compact sets and strip type sets; monotonicity of part metric between
two positive solutions of \eqref{main-eq};  the existence, uniqueness, and stability of time and space periodic positive solutions of \eqref{periodic-eq}
and spatial spreading speeds of \eqref{periodic-eq}; and the principal eigenvalues theory for time periodic dispersal operators.

\subsection{Comparison principle and global existence}

In this subsection, we consider comparison principle and global existence of solutions of \eqref{main-eq}.
Throughout this subsection, we assume (H0).

Let $\Omega\subset\mathcal{H}$ be a convex region of $\mathcal{H}$.
For a given continuous and bounded function $u:[0,\tau)\times\RR^N\to\RR$, it is called a {\it super-solution} ({\it sub-solution}) of \eqref{main-eq} on $[0,\tau)\times\Omega$ if
\begin{equation}
\label{sub-super-solution-eq}
u_t(t,x)\geq (\leq) (\mathcal{A} u)(t,x)+u(t,x)f(t, x,u(t,x))\quad\forall  (t,x)\in (0,\tau)\times \bar\Omega.
\end{equation}

\begin{proposition}[Comparison principle]
\label{basic-comparison}
\begin{itemize}
\item[(1)]  Suppose that $u^1(t,x)$ and $u^2(t,x)$ are sub- and super-solutions of \eqref{main-eq}
on $[0,\tau)\times\Omega$ with $u^1(t,x)\leq u^2(t,x)$ for $x\in \mathcal{H}\setminus \Omega$, $t\in[0,\tau)$ and $u^1(0,x)\leq u^2(0,x)$ for $x\in \bar \Omega$. Then
$u^1(t,x)\leq u^2(t,x)$ for  $x\in \Omega$ and $t\in [0,\tau)$. Moreover, if $u^1(0,x)\not\equiv u^2(0,x)$ for $x\in\Omega$, then
$u^1(t,x)<u^2(t,x)$ for $t\in(0,\tau)$ and $x\in\Omega$.

\item[(2)] If $u_{01},u_{02}\in X$ and $u_{01}\leq u_{02}$, then $u(t,\cdot;u_{01})\leq u(t,\cdot;u_{02})$ for $t>0$ at which both
$u(t,\cdot;u_{01})$ and $u(t,\cdot;u_{02})$  exist. Moreover, if  $u_{01}\not =u_{02}$, then $u(t,x;u_{01})< u(t,x;u_{02})$ for all $x\in\mathcal{H}$
 and $t>0$ at which both
$u(t,\cdot;u_{01})$ and $u(t,\cdot;u_{02})$  exist.

\item[(3)] If $u_{01},u_{02}\in X$ and $u_{01}\ll u_{02}$, then  $u(t,\cdot;u_{01})\ll u(t,\cdot;u_{02})$ for $t>0$ at which both
$u(t,\cdot;u_{01})$ and $u(t,\cdot;u_{02})$  exist.
\end{itemize}
\end{proposition}

\begin{proof}
(1) The case that $\mathcal{H}=\mathcal{H}_1(=\RR^N)$ and $\mathcal{A}=\mathcal{H}_1(=\Delta)$ follows from comparison principle for parabolic equations.
We prove the case that $\mathcal{H}=\mathcal{H}_2(=\RR^N)$ and $\mathcal{A}=\mathcal{A}_2$. The case
that $\mathcal{H}=\mathcal{H}_3(=\ZZ^N)$ and $\mathcal{A}=\mathcal{A}_3$  can be proved similarly.

Observe that for any $t\in[0,\tau)$,
\begin{align}
\label{comparison-eq1}
\int_{\RR^N}\kappa(y-x)u^1(t,y)dy&=\int_{\RR^N\setminus\Omega}\kappa(y-x)u^1(t,y)dy+\int_{\Omega}\kappa(y-x)u^1(t,y)dy\nonumber\\
&\le\int_{\RR^N\setminus\Omega}\kappa(y-x)u^2(t,y)dy+\int_\Omega\kappa(y-x)u^1(t,y)dy.
\end{align}
Let $v(t,x)=u^2(t,x)-u^1(t,x)$. By \eqref{comparison-eq1},
\begin{align*}
v_t(t,x)&\ge \int_{\Omega}\kappa(y-x)v(t,y)dy-v(t,x)+u^2(t,x)f(t,x,u^2(t,x))-u^1(t,x)f(t,x,u^1(t,x))\\
&=\int_{\Omega}\kappa(y-x)v(t,y)dy-v(t,x)+a(t,x)v(t,x),\quad x\in\bar \Omega,\,\, t\in(0,\tau),
\end{align*}
where
\begin{equation*}
a(t,x)=f(t,x,u^2(t,x))+u^1(t,x)\int_0^1 \p_u f(t,x,s u^2(t,x)+(1-s)u^1(t,x))ds.
\end{equation*}
The rest of the proof follows from the arguments of \cite[Proposition 2.4]{HuShVi}.

(2) It follows from (1) with $u^1(t,x)=u(t,x;u_{01})$, $u^2(t,x)=u(t,x;u_{02})$, and $\Omega=\mathcal{H}$.

(3) We provide a proof for the case that $\mathcal{H}=\mathcal{H}_2$ and $\mathcal{A}=\mathcal{A}_2$.
Other cases can be proved similarly.
Take any $\tau>0$ such that both $u(t,\cdot;u_{01})$ and $u(t,\cdot;u_{02})$ exist on
$[0,\tau]$. It suffices to prove that $u(t,\cdot;u_{02})\gg u(t,\cdot;u_{01})$
for $t\in [0,\tau]$. To this end, let $w(t,x)=u(t,x;u_{02})-u(t,x;u_{01})$. Then $w(t,x)$ satisfies the following equation,
\begin{equation*}
w_t(t,x)=\int_{\RR^N}\kappa(y-x)w(t,y)dy-w(t,x)+a(t,x)w(t,x),
\end{equation*}
where
\begin{align*}
a(t,x)=f(t,x,u(t,x;u_{02}))+u(t,x;u_{01})\int_0^1 \p_u f(t,x,s u(t,x;u_{02})+(1-s)u(t,x;u_{01}))ds.
\end{align*}
Let $M>0$ be such that
$M\geq \sup_{x\in\RR^N,t\in[0,\tau]}(1-a(t,x))$
and $\tilde w(t,x)=e^{Mt}w(t,x)$.
Then $\tilde w(t,x)$ satisfies
$$
\tilde w_t(t,x)=\int_{\RR^N}\kappa(y-x)\tilde w(t,y)dy+[M-1+a(t,x)]\tilde w(t,x).
$$
Let $\mathcal{K}:X\to X$ be defined by
\begin{equation}
\label{k-op}
(\mathcal{K}u)(x)=\int_{\RR^N}\kappa(y-x)u(y)dy\quad {\rm for}\quad u\in X.
\end{equation}
Then $\mathcal{K}$ generates an analytic semigroup on $X$ and
$$
\tilde w(t,\cdot)=e^{\mathcal{K}t}(u_{02}-u_{01})+\int_0^t e^{\mathcal{K}(t-\tau)}(M-1+a(\tau,\cdot))\tilde w(\tau,\cdot)d\tau.
$$
Observe that $e^{\mathcal{K} t}u_0\ge 0$ for any $u_0\in X^+$ and $t\ge 0$ and
$e^{\mathcal{K}t}u_0\gg 0$ for any $u_0\in X^{++}$ and $t\ge 0$. Observe also that
$u_{02}-u_{01}\in X_2^{++}$. By (2), $\tilde w(\tau,\cdot)\geq 0$ and hence $(M-1+a(\tau,\cdot))\tilde w(\tau,\cdot)\ge 0$
for $\tau\in[0,T]$.  It then follows that
$\tilde w(t,\cdot)\gg 0$ and then $w(t,\cdot)\gg 0$ (i.e. $u(t,\cdot;u_{02})\gg u(t,\cdot;u_{01})$) for $t\in [0,\tau]$.
\end{proof}

\begin{proposition}[Global existence]
\label{basic-global-existence}
For any given  $u_0\in X^+$, $u(t,\cdot;u_0)$ exists for all $t\geq 0$.
\end{proposition}

\begin{proof}
Let  $u_0\in X^+$ be given. There is $M\gg 1$ such that
$0\leq u_0(x)\leq M$ and $f(t,x,M)<0$ for all $x\in\mathcal{H}$. Then by Proposition \ref{basic-comparison},
$$
0\leq u(t,\cdot;u_0)\leq M
$$
for any $t>0$ at which $u(t,\cdot;u_0)$ exists. It is then not difficult to prove that
for any $\tau>0$ such that $u(t,\cdot;u_0)$ exists on $(0,\tau)$, $\lim_{t\to \tau}u(t,\cdot;u_0)$ exists
in $X$.
 This implies that $u(t,\cdot;u_0)$ exists and
$u(t,\cdot;u_0)\geq 0$ for all $t\geq 0$.
\end{proof}

\subsection{Convergence on compact subsets and strip type subsets}

In this subsection, we explore the convergence property of solutions of \eqref{main-eq} on compact subsets and strip type subsets.
As mentioned before, to indicate the dependence of solutions of \eqref{main-eq} on the nonlinearity,
we may write $u(t,\cdot;u_0)$ as $u(t,\cdot;u_0,f)$.

\begin{proposition}[Convergence on compact and strip type subsets]
\label{basic-convergence}
Suppose that $u_{0n},u_0\in X^+$ $(n=1,2,\cdots)$ with $\{\|u_{0n}\|\}$ being  bounded, and $f_n$, $g_n$ ($n=1,2 \cdots$) satisfy (H0) with
$f_n(t,x,u)$, $g_n(t,x,u)$, and $\p_uf_n(t,x,u)$ being bounded uniformly in $x\in\mathcal{H}$ and $(t,u)$ on bounded subsets.
\begin{itemize}
\item[(1)]
If $u_{0n}(x)\to u_0(x)$ as $n\to\infty$
 uniformly in $x$ on bounded sets and $f_n(t,x,u) - g_n(t,x,u)$ as $n\to\infty$ uniformly  in $(t,x,u)$ on bounded sets,
 then for each $t>0$, $u(t,x;u_{0n},f_n)- u(t,x;u_0,g_n)\to 0$ as $n\to\infty$ uniformly in $x$ on bounded sets.

\item[(2)]
If $u_{0n}(x)\to u_0(x)$ as $n\to\infty$
 uniformly in $x$ on  any set $E$ with $\{x\cdot\xi\,|\, x\in E\}$ being a bounded set of $\RR$  and $f_n(t,x,u) - g_n(t,x,u)\to 0$ as $n\to\infty$ uniformly  in $(t,x,u)$ on any set $E$ with
 $\{(t,x\cdot\xi,u)\,|\, (t,x,u)\in E\}$ being a bounded  set of $\RR^3$,
 then for each $t>0$, $u(t,x;u_{0n},f_n)-u(t,x;u_0,g_n)\to 0$ as $n\to\infty$ uniformly in $x$  on any set $E$ with $\{x\cdot\xi\,|\, x\in E\}$ being a bounded set of $\RR$.
\end{itemize}
\end{proposition}

\begin{proof}
(1)
We prove the case that $\mathcal{H}=\mathcal{H}_2$ and $\mathcal{A}_2$. Other cases can be proved similarly.

 Let $v^n(t,x)=u(t,x;u_{0n},f_n)-u(t,x;u_0,g_n)$.
Then $v^n(t,x)$ satisfies
\begin{equation*}
v^n_t(t,x)=\int_{\RR^N}\kappa(y-x)v^n(t,y)dy-v^n(t,x)+a_n(t,x)v^n(t,x)+b_n(t,x),
\end{equation*}
where
\begin{align*}
a_n(t,x)=&f_n(t,x,u(t,x;u_{0n},f_n))\\
&+u(t,x;u_{0},g_n)\cdot \int_0^1 \p_u f_n(t,x,s u(t,x;u_{0n},f_n)+
(1-s)u(t,x;u_0,g_n))ds
\end{align*}
and
\begin{align*}
b_n(t,x)= u(t,x;u_0,g_n)
\cdot \big(f_n(t,x, u(t,x;u_0,g_n))-g_n(t,x,
 u(t,x;u_0,g_n))\big).
\end{align*}
Observe that $\{a_n(t,x)\}$ is uniformly bounded and continuous  in $t$ and $x$ and
$b_n(t,x)\to 0$ as $n\to\infty$ uniformly in $(t,x)$ on bounded sets of $[0,\infty)\times\RR^N$.

Take a $\rho>0$. Let
$$
X(\rho)=\{u\in C(\RR^N,\RR)\,|\, u(\cdot) e^{-\rho \|\cdot\|}\in X\}
$$
with norm $\|u\|_\rho=\|u(\cdot)e^{-\rho \|\cdot\|}\|$.
Note that $\mathcal{K}:X(\rho)\to X(\rho)$ also generates an analytic  semigroup,
where $\mathcal{K}$ is as in \eqref{k-op},
and there are $M>0$ and $\omega>0$ such that
$$
\|e^{(\mathcal{K}-\mathcal{I})t}\|_{X(\rho)}\leq M e^{\omega t}\quad \forall t\geq 0,
$$
where $\mathcal{I}$ is the identity map on $X(\rho)$.
Hence
\begin{align*}
v^n(t,\cdot)=&e^{(\mathcal{K}-\mathcal{I})t}v^n(0,\cdot)+\int_0^t e^{(\mathcal{K}-\mathcal{I})(t-\tau)}a_n(\tau,\cdot)v^n(\tau,\cdot)
d\tau\\
&+\int_0^t e^{(\mathcal{K}-\mathcal{I})(t-\tau)}b_n(\tau,\cdot)d\tau
\end{align*}
and
then
\begin{align*}
\|v^n(t,\cdot)\|_{X(\rho)}&\leq M e^{\omega t}\|v^n(0,\cdot)\|_{X(\rho)}+
M\sup_{\tau\in[0,t],x\in\RR^N} |a_n(\tau,x)|\int_0^ t e^{\omega(t-\tau)}\|v^n(\tau,\cdot)\|_{X(\rho)}d\tau\\
&\quad +M\int_0^ t e^{\omega(t-\tau)}\|b_n(\tau,\cdot)\|_{X(\rho)}d\tau\\
&\leq M e^{\omega t}\|v^n(0,\cdot)\|_{X(\rho)}+
M\sup_{\tau\in[0,t],x\in\RR^N} |a_n(\tau,x)|\int_0^ t e^{\omega(t-\tau)}\|v^n(\tau,\cdot)\|_{X(\rho)}d\tau\\
&\quad + \frac{M}{\omega}\sup_{\tau\in
[0,t]}\|b_n(\tau,\cdot)\|_{X(\rho)} e^{\omega t}.
\end{align*}
By Gronwall's inequality,
$$
\|v^n(t,\cdot)\|_{X(\rho)}\leq e^{(\omega+M\sup_{\tau\in[0,t],x\in\RR^N} |a_n(\tau,x)|)t}\Big(M\|v^n(0,\cdot)\|_{X(\rho)}
+ \frac{M}{\omega}\sup_{\tau\in [0,t]}\|b_n(\tau,\cdot)\|_{X(\rho)}\Big).
$$
Note that $\|v^n(0,\cdot)\|_{X(\rho)}\to 0$ and $\sup_{\tau\in
[0,t]}\|b_n(\tau,\cdot)\|_{X(\rho)}\to 0$ as $n\to\infty$. It then
follows that
$$
\|v^n(t,\cdot)\|_{X(\rho)}\to 0\quad {\rm as}\quad n\to\infty
$$ and then
$$
u(t,x;u_{0n},f_n)-u(t,x;u_0,g_n)\quad {\rm as}\quad
n\to\infty
$$
uniformly in $x$ on bounded sets.

(2)  It can be proved by similar arguments as in (1) with $X(\rho)$ being replaced by
$X_{\xi}(\rho)$, where
$$
X_{\xi}(\rho)=\{u\in C(\mathcal{H},\RR)\,|\, u_{\xi,\rho}\in X\},
$$
with norm $\|u\|_{X_{\xi}(\rho)}=\|u_{\xi,\rho}\|_X$,
where $u_{\xi,\rho}(x)=e^{-\rho |x\cdot\xi|}u(x)$.
\end{proof}

\subsection{Part metric}

In this subsection, we investigate the decreasing property of the so called part metric between two positive solutions of
\eqref{main-eq}. Throughout this subsection, we also assume (H0).

First, we introduce the notion of part metric.
For given $u,v\in X^{++}$, define
$$
\rho(u,v)=\inf\{\ln \alpha\,|\, \frac{1}{\alpha}u\leq v\leq \alpha u,\,\, \alpha\geq 1\}.
$$
Observe that $\rho(u,v)$ is well defined and there is $\alpha\geq 1$ such that $\rho(u,v)=\ln\alpha$.
Moreover, $\rho(u,v)=\rho(v,u)$ and $\rho(u,v)=0$ iff $u\equiv v$. In literature, $\rho(u,v)$ is called the {\it part metric}
between $u$ and $v$.

We remark that the concept of part metric was introduced by Thompson  \cite{Tho}. It was first observed by Krause and Nussbaum \cite{KrNu}
that a monotone map in a strong ordered Banach space  with strong subhomogeneity is contractive with respect to the part metric on the interior of the positive cone
of the Banach space.

\begin{proposition}[Strict decreasing of part metric]
\label{basic-part-metric}
For any $\epsilon>0$, $\sigma>0$, $M>0$,  and $\tau>0$  with $\epsilon<M$ and
$\sigma\le \ln \frac{M}{\epsilon}$, there is $\delta>0$   such that
for any $u_0,v_0\in X^{++}$ with $\epsilon\le u_0(x)\le M$, $\epsilon\le v_0(x)\le M$ for $x\in\mathcal{H}$ and
$\rho(u_0,v_0)\ge\sigma$, there holds
$$
\rho(u(\tau,\cdot;u_0),u(\tau,\cdot;v_0))\le \rho(u_0,v_0)- \delta.
$$
\end{proposition}

\begin{proof}
We give a proof for the case that $\mathcal{H}=\mathcal{H}_1$ and $\mathcal{A}=\mathcal{A}_1$. Other cases can be proved similarly.

Let $\epsilon>0$, $\sigma>0$, $M>0$, and $\tau>0$ be given and $\epsilon<M$, $\sigma<\ln \frac{M}{\epsilon}$. First, note that by Proposition \ref{basic-comparison},
 there are $\epsilon_1>0$ and $M_1>0$ such that
for any $u_0\in X^{++}$ with $\epsilon\le u_0(x)\le M$ for $x\in\RR^N$, there holds
\begin{equation}
\label{part-metric-eq1}
\epsilon_1\le u(t,x;u_0)\le M_1\quad \forall\,\, t\in[0,\tau],\,\, x\in\RR^N.
\end{equation}
Let
\begin{equation}
\label{part-metric-eq2}
\delta_1=\epsilon_1^2 e^\sigma (1-e^\sigma)\sup_{t\in[0,\tau],x\in\RR^N,u\in[\epsilon_1,M_1M/\epsilon]}f_u(t,x,u).
\end{equation}
Then $\delta_1>0$ and there is $0<\tau_1\le\tau$ such that
\begin{equation}
\label{part-metric-eq3-1}
\frac{\delta_1}{2}\tau_1<e^\sigma\epsilon_1
\end{equation}
and
\begin{equation}
\label{part-metric-eq3-2}
\Big|\frac{\delta_1}{2}tvf_u(t,x,w)\Big|+\Big|\frac{\delta_1}{2}tf(t,x,v-\frac{\delta_1}{2}t)\Big|\le\frac{\delta_1}{2}\quad \forall \,\, t\in [0,\tau_1],\,\, x\in\RR^N,\,\, v,w\in[0,M_1M/\epsilon].
\end{equation}
Let
\begin{equation}
\label{part-metric-eq4}
\delta_2=\frac{\delta_1\tau_1}{2 M_1}.
\end{equation}
Then $\delta_2<e^\sigma$ and $0<\frac{\delta_2 \epsilon}{M}<1$. Let
\begin{equation}
\label{part-metric-eq5}
\delta=-\ln\big(1-\frac{\delta_2\epsilon}{M}\big).
\end{equation}
Then $\delta>0$. We prove that  $\delta$ defined in \eqref{part-metric-eq5}  satisfies the property in the proposition.

For any $u_0,v_0\in X^{++}$ with $\epsilon\le u_0(x)\le M$ and $\epsilon\le v_0(x)\le M$ for $x\in\RR$ and
$\rho(u_0,v_0)\ge\sigma$,  there is $\alpha^*\ge 1$ such that  $\rho(u_0,v_0)=\ln \alpha^{*}$
and $\frac{1}{\alpha^{*}}u_0\leq v_0\leq \alpha^{*} u_0$. Note that $e^\sigma\le\alpha^*\le \frac{M}{\epsilon}$.
We first show that $\rho(u(t,\cdot;u_0),u(t,\cdot;v_0))$ is non-increasing in $t>0$.

By Proposition \ref{basic-comparison},
$$u(t, \cdot;v_0)\le u(t, \cdot; \alpha^{*}u_0)\quad {\rm for}\quad t>0.
$$
Let
 $$v(t, x)=\alpha^{*}u(t, x; u_0).$$
We then have
\begin{align*}
v_{t}(t,x)&=\Delta v(t,x)+v(t,x)f(t,x, u(t,x;u_0))\\
& =\Delta
v(t,x)+v(t,x)f(t, x,v(t,x))+v(t,x)f(t, x,u(t,x;u_0))-v(t,x)f(t,x,v(t,x))\\
&\ge \Delta v (t,x)+ v(t,x)f(t, x,v(t,x))\quad \forall t>0.
\end{align*}  By Proposition \ref{basic-comparison} again,
$$
u(t,\cdot;\alpha^* u_0)\leq \alpha^* u(t,\cdot;u_0)
$$
and hence
$$
u(t,\cdot;v_0)\leq \alpha^* u(t,\cdot;u_0)
$$
for $t>0$.
Similarly, we can prove that
$$
\frac{1}{\alpha^*}u(t,\cdot;u_0)\le u(t,\cdot;v_0)
$$
for $t>0$. It then follows that
$$
\rho(u(t,\cdot;u_0),u(t,\cdot;v_0))\le \rho(u_0,v_0)\quad \forall \,\, t\ge 0
$$
and then
$$
\rho(u(t_2,\cdot;u_0),u(t_2,\cdot;v_0))\le \rho(u(t_1,\cdot;u_0),u(t_1,\cdot;v_0))\quad \forall\,\, 0\le t_1\le t_2.
$$

Next, we prove that
$$
\rho(u(\tau,\cdot;u_0),u(\tau,\cdot;v_0))\le\rho(u_0,v_0)-\delta.
$$
Note that $e^\sigma\le \alpha^*\le \frac{M}{\epsilon}$
and
\begin{align*}
v_{t}(t,x)&=\Delta v(t,x)+v(t,x)f(t,x, u(t,x;u_0))\\
& =\Delta
v(t,x)+v(t,x)f(t, x,v(t,x))+v(t,x)f(t, x,u(t,x;u_0))-v(t,x)f(t,x,v(t,x))\\
&\ge \Delta v (t,x)+ v(t,x)f(t, x,v(t,x))+\delta_1\quad \forall 0<t\le\tau_1.
\end{align*}  This together with \eqref{part-metric-eq3-1}, \eqref{part-metric-eq3-2}
implies that
$$
(v(t,x)-\frac{\delta_1}{2}t)_t\ge \Delta\big(v(t,x)-\frac{\delta_1}{2}t\big)+\big(v(t,x)-\frac{\delta_1}{2}t\big)f\big(t,x,v(t,x)-\frac{\delta_1}{2}t\big)
$$
for $0<t\le \tau_1$.
Then by Proposition \ref{basic-comparison} again,
$$
u(t,\cdot;\alpha^* u_0)\leq \alpha^* u(t,\cdot;u_0)-\frac{\delta_1}{2}t\quad {\rm for}\quad 0<t\le\tau_1.
$$
By \eqref{part-metric-eq4},
$$
u(\tau_1,\cdot;v_0)\le (\alpha^*-\delta_2) u(\tau_1,\cdot;u_0).
$$
Similarly, it can be proved that
$$
\frac{1}{\alpha^*-\delta_2}u(\tau_1,\cdot;u_0)\le u(\tau_1,\cdot;v_0).
$$
It then follows that
$$
\rho(u(\tau_1,\cdot;u_0),u(\tau_1,\cdot;v_0))\le \ln(\alpha^*-\delta_2) =\ln\alpha^*+\ln (1-\frac{\delta_2}{\alpha^*})\le \rho(u_0,v_0)-\delta.
$$
and hence
$$
\rho(u(\tau,\cdot;u_0),u(\tau,\cdot;v_0))\le\rho(u(\tau_1,\cdot;u_0),u(\tau_1,\cdot;v_0))\le\rho(u_0,v_0)-\delta.
$$
\end{proof}

\subsection{Time periodic positive solutions and spreading speeds of KPP equations in periodic media}

In this subsection, we recall some existing results on the existence, uniqueness, and stability of time and space periodic
positive solutions and spatial spreading speeds of \eqref{periodic-eq}.

A solution $u(t,x)$ of \eqref{periodic-eq}   is called {\it time and space periodic solution} if it is a solution on $t\in\RR$ and  $u(t+T,x)=u(t,x+p_i{\bf e_i})=u(t,x)$ for
$t\in\RR$, $x\in\mathcal{H}$, and $i=1,2,\cdots,N$. It is called a {\it positive solution} if $u(t,x)>0$ for all $t$ in the existence interval and $x\in\mathcal{H}$.

\begin{proposition}
\label{basic-periodic-solution}
Consider \eqref{periodic-eq} and assume that $f_0$ satisfies (H0) and $f_0(t,x+p_i{\bf e_i},u)=f_0(t,x,u)$ ($i=1,2,\cdots,N$).
\begin{itemize}
\item[(1)] (Uniqueness of periodic positive solutions) If \eqref{periodic-eq} has a time and space periodic positive solution, then it is unique.

\item[(2)] (Stability of periodic positive solutions) If \eqref{periodic-eq} has a time and space periodic positive solution $u^*(t,x)$, then
it is globally asymptotically stable with respect to perturbations in $X^{+}_p\setminus\{0\}$.

\item[(3)] (Existence of periodic positive solutions) If $\lambda(f_0(\cdot,\cdot,0))>0$, then \eqref{periodic-eq} has a time and space periodic positive
solution.
\end{itemize}
\end{proposition}

\begin{proof}
The case that $\mathcal{H}=\mathcal{H}_1$ and $\mathcal{A}=\mathcal{A}_1$ follows from  Theorems 1.1, 1.3, and 1.6 in \cite{Nad0}.
The case that $\mathcal{H}=\mathcal{H}_2$ and $\mathcal{A}=\mathcal{A}_2$ follows from Theorem E in \cite{RaSh}.
The case that $\mathcal{H}=\mathcal{H}_3$ and $\mathcal{A}=\mathcal{A}_3$  can be proved by the similar
arguments as in Theorem E in \cite{RaSh}.
\end{proof}

We remark that if $f_0$ satisfies (H0), then $g_0(t,x,u):=uf_0(t,x,u)$ is strictly subhomogeneous (or sublinear in some literature) in $u$ in the sense that
$g_0(t,x,su)>g(t,x,u)$ for all $t\in\RR$, $x\in\mathcal{H}$, $u>0$ and $s\in (0,1)$. The reader is referred to \cite[Section 2.3]{Zha2} for general theories on subhomogeneous systems.

\begin{proposition}
\label{basic-spreading-speed} Consider \eqref{periodic-eq}. Assume
that $f_0$ satisfies (H0), $f_0(t,x+p_i{\bf e_i},u)=f_0(t,x,u)$ ($i=1,2,\cdots,N$),  and $\lambda(f_0(\cdot,\cdot,0))>0$.
Then for any given $\xi\in S^{N-1}$, \eqref{periodic-eq} has a spatial spreading speed $c_0^*(\xi)$ in the direction of $\xi$.
Moreover,
$$
c_0^*(\xi)=\inf_{\mu>0}\frac{\lambda_{\xi,\mu}(f_0(\cdot,\cdot,0))}{\mu}
$$
and for any $c<c_0^*(\xi)$ and $u_0\in X^+(\xi)$,
$$
\limsup_{x\cdot\xi\le ct,t\to\infty}|u(t,x;u_0,f_0(\cdot,\cdot+z,\cdot))-u_0^*(t,x+z)|=0
$$
uniformly in $z\in\mathcal{H}$.
\end{proposition}

\begin{proof}
The cases that $\mathcal{H}=\mathcal{H}_i$ and $\mathcal{A}=\mathcal{A}_i$ for $i=1,3$ follow
from Theorem 2.1  in \cite{Wei2} (see also \cite{LiZh2}, \cite{NoRuXi}, \cite{NoXi}).
The case that $\mathcal{H}=\mathcal{H}_2$ and $\mathcal{A}=\mathcal{A}_2$ follows from Theorem 4.1 in \cite{RaShZh}.
\end{proof}

\begin{proposition}
\label{basic-spreading-feature}
Consider \eqref{periodic-eq} and assume that the conditions in Proposition \ref{basic-spreading-speed} hold. Then
for any given $\xi\in S^{N-1}$, the following hold.
\begin{itemize}
\item[(1)] For each $ u_0\in X^+$ satisfying that  $u_0(x)=0$ for $x\in\mathcal{H}$ with $|x\cdot\xi|\gg 1$,
$$
\limsup_{|x\cdot\xi|\geq ct, t\to\infty}u(t,x;u_0,f_0(\cdot,\cdot+z,\cdot))=0 \quad \forall\,\, c>\max\{c_0^*(\xi),c_0^*(-\xi)\}
$$
uniformly in $z\in\mathcal{H}$.

\item[(2)]  For each $\sigma>0$,  $r>0$, and  $u_0\in X^+$ satisfying that $u_0(x)\geq\sigma$ for $x\in\mathcal{H}$ with $|x\cdot \xi|\leq r$,
$$
\limsup_{|x\cdot \xi|\leq ct,t\to\infty}|u(t,x;u_0,f_0(\cdot,\cdot+z,\cdot))-u^*(t, x+z)|=0\quad \forall\,\, 0<c<\min\{c_0^*(\xi),c_0^*(-\xi)\}
$$
uniformly in $z\in\mathcal{H}$.

\item[(3)] For each $u_0\in X^+$ satisfying that $u_0(x)=0$ for $x\in\mathcal{H}$ with $\|x\|\gg 1$,
$$
\limsup_{\|x\|\geq ct,t\to\infty}u(t,x;u_0,f_0(\cdot,\cdot+z,\cdot))=0\quad \forall\,\, c>\sup_{\xi\in S^{N-1}}c_0^*(\xi)
$$
uniformly in $z\in\mathcal{H}$.

\item[(4)] For each $\sigma>0$, $r>0$, and  $u_0\in X^+$ satisfying that $u_0(x)\geq \sigma$ for $\|x\|\le r$,
$$
\limsup_{\|x\|\leq ct,t\to\infty}|u(t,x;u_0,f_0(\cdot,\cdot+z,\cdot))-u_0^*(t, x+z)|=0\quad \forall\,\, 0<c<\inf_{\xi\in S^{N-1}}c_0^*(\xi)
$$
uniformly in $z\in\mathcal{H}$.
\end{itemize}
\end{proposition}

\begin{proof}
The cases that $\mathcal{H}=\mathcal{H}_i$ and $\mathcal{A}=\mathcal{A}_i$ for $i=1,3$ follow
from the arguments of Theorems 2.1-2.3 in \cite{Wei2} (see also \cite{LiZh2}),
and the case that $\mathcal{H}=\mathcal{H}_2$ and $\mathcal{A}=\mathcal{A}_2$ follows from the arguments of Theorem 4.2  in \cite{RaShZh}.
\end{proof}
\subsection{Principal eigenvalues of time periodic dispersal operators}

In this subsection, we recall some principal eigenvalues theory for time periodic dispersal operators.

Let $\mathcal{X}_p$ and $\mathcal{A}_{\xi,\mu}$ be as in \eqref{lift-x-space-eq} and \eqref{operator-a-xi-mu-eq}, respectively
($\xi\in S^{N-1}$ and $\mu\in\RR$). For given $a\in \mathcal{X}_p$, $\xi\in S^{N-1}$,  and $\mu\in\RR$, let $\lambda_{\xi,\mu}(a)$
be as in \eqref{lambda-max-eq}.

\begin{definition}
A real number $\lambda_0$ is said to be the {\rm principal eigenvalue} of $-\p_t+\mathcal{A}_{\xi,\mu}+a(\cdot,\cdot)\mathcal{I}$
if $\lambda_0$ is an isolated eigenvalue of  $-\p_t+\mathcal{A}_{\xi,\mu}+a(\cdot,\cdot)\mathcal{I}$ with
 a positive eigenfunction $\phi(\cdot,\cdot)$ (i.e. $\phi(t,x)>0$ for $(t,x)\in\RR\times\mathcal{H}$ and $\phi(\cdot,\cdot)\in\mathcal{X}_p$)
and for any $\lambda\in \sigma(-\p_t+\mathcal{A}_{\xi,\mu}+a(\cdot,\cdot)\mathcal{I})$, ${\rm Re}\mu\le \lambda_0$.
\end{definition}

We remark that $\lambda_{\xi,\mu}(a)\in \sigma(-\p_t+\mathcal{A}_{\xi,\mu}+a(\cdot,\cdot)\mathcal{I})$ and if $-\p_t+\mathcal{A}_{\xi,\mu}+a(\cdot,\cdot)\mathcal{I}$
admits a principal eigenvalue $\lambda_0$, then $\lambda_0=\lambda_{\xi,\mu}(a)$.
We also remark that in the case that $\mathcal{H}=\mathcal{H}_i$ and $\mathcal{A}=\mathcal{A}_i$ with $i=1$ or $3$, principal eigenvalue
of $-\p_t+\mathcal{A}_{\xi,\mu}+a(\cdot,\cdot)\mathcal{I}$ always exists. But in the case that $\mathcal{H}=\mathcal{H}_2$ and $\mathcal{A}=\mathcal{A}_2$,
$-\p_t+\mathcal{A}_{\xi,\mu}+a(\cdot,\cdot)\mathcal{I}$ may not have a principal eigenvalue (see \cite{Cov} and \cite{ShZh1} for examples).
The following proposition is established in \cite{RaSh} regarding the existence of  principal eigenvalues of time periodic nonlocal dispersal
operators.

For given $a\in\mathcal{X}_p$, let
$$
\hat a(x)=\frac{1}{T}\int_0^T a(t,x)dt.
$$

\begin{proposition}
\label{basic-criterion}
\begin{itemize}
\item[(1)]
If $\hat a(\cdot)$ is $C^N$ and there is $x_0\in\RR^N$ such that $\hat a(x_0)=\max_{x\in\RR^N}\hat a(x_0)$ and the partial derivatives
of $\hat a(x)$ up to order $N-1$ at $x_0$ are zero, then for any $\xi\in S^{N-1}$ and $\mu\in\RR$, $\lambda_{\xi,\mu}(a)$ is
the principal eigenvalue of  $-\p_t+\mathcal{A}_{\xi,\mu}+a(\cdot,\cdot)\mathcal{I}$.

\item[(2)] Let $a(\cdot,\cdot)\in \mathcal{X}_p$ be given. For any $\epsilon>0$, there is $a^\pm (\cdot,\cdot)\in\mathcal{X}_p$ such that
$\lambda_{\xi,\mu}(a^\pm)$ are principal eigenvalues of  $-\p_t+\mathcal{A}_{\xi,\mu}+a^\pm(\cdot,\cdot)\mathcal{I}$,
$$
a^-(t,x)\le a(t,x)\le a^+(t,x)\quad \forall (t,x)\in\RR\times\mathcal{H},
$$
and
$$
\lambda_{\xi,\mu}(a^+)-\epsilon\le\lambda_{\xi,\mu}(a)\le \lambda_{\xi,\mu}(a^-)+\epsilon.
$$
\end{itemize}
\end{proposition}

\begin{proof}
We only need to prove the case that $\mathcal{H}=\mathcal{H}_2$ and $\mathcal{A}=\mathcal{A}_2$.

(1) It follows from \cite[Theorem B(1)]{RaSh}.

(2) It follows from \cite[Proposition 3.10, Lemma 4.1]{RaSh}.
\end{proof}

\section{Existence, Uniqueness, and Stability of Time Periodic Strictly Positive Solutions}

In this section, we explore the existence, uniqueness, and stability of time periodic strictly
positive solutions of \eqref{main-eq} and prove Theorem \ref{positive-solution-thm}.

\subsection{Uniqueness and stability}

In this subsection, we prove the uniqueness and stability of time periodic strictly positive solutions of
\eqref{main-eq} (if exist), i.e. prove Theorem \ref{positive-solution-thm}(1) and (2).

\begin{proof}[Proof of Theorem \ref{positive-solution-thm} (1)]

Suppose that there are two time periodic strictly positive solutions $u^1(t,x)$ and
$u^2(t,x)$.
Then
$$
\rho(u^1(t+T,\cdot),u^2(t+T,\cdot))=\rho(u^1(t,\cdot),u^2(t,\cdot))
$$
for any $t\in\RR$. By Proposition \ref{basic-part-metric}, we must have
$$u^1(t,x)\equiv u^2(t,x).$$
Thus time periodic positive solutions are unique.
\end{proof}

\begin{proof} [Proof of Theorem \ref{positive-solution-thm}  (2)]

First of all, for any $u_0\in X^{++}$, by Proposition \ref{basic-part-metric},
\begin{equation}
\label{stability-eq1}
\rho(u(t,\cdot;u_0),u^*(t,\cdot))\le \rho(u_0,u^*(0,\cdot))\quad \forall \,\, t\ge 0.
\end{equation}
This implies that $u^*(t,x)$ is stable with respect to perturbations in $X^{++}$.

Next, for any given $u_0\in X^{++}$, we show
\begin{equation}
\label{stability-aux-eq1}
\|u(t,\cdot;u_0,f(\cdot+\tau,\cdot,\cdot))-u^*(t+\tau,\cdot)\|\to 0
\end{equation}
as $t\to\infty$ uniformly in $\tau\in\RR$. Thanks to the periodicity
of $f(t,x,u)$ and $u^*(t,x)$ in $t$, we only need to show that the
limit in \eqref{stability-aux-eq1} exists and is uniformly in
$\tau\in[0,T]$. Moreover, note that
$$
u(t,x;u_0,f(\cdot+\tau,\cdot,\cdot))=u(t-T+\tau,x;u(T-\tau,\cdot;u_0,f(\cdot+\tau,\cdot,\cdot)),f)\quad
\forall t\ge T-\tau.
$$
By Proposition \ref{basic-comparison},
$$
\inf_{\tau\in[0,T],x\in\mathcal{H}}u(T-\tau,x;u_0,f(\cdot+\tau,\cdot,\cdot))>0.
$$
It then suffices to prove that the limit in \eqref{stability-aux-eq1}
exists for $\tau=0$.

Let $\alpha_0\ge 1$ be such that $\rho(u_0,u^*(0,\cdot))=\ln
\alpha_0$ and
$$
\frac{1}{\alpha_0}u^*(0,x)\le u_0(x)\le \alpha_0 u^*(0,x)\quad \forall \,\, x\in\mathcal{H}.
$$
By Proposition \ref{basic-part-metric},
there is $\alpha_\infty\ge 1$ such that
$$
\lim_{t\to\infty}\rho(u(t,\cdot;u_0),u^*(t,\cdot))=\ln\alpha_\infty.
$$
Moreover, by \eqref{stability-eq1},
$$
\rho(u(t,\cdot;u_0),u^*(t,\cdot))\le \rho(u_0,u^*(0,\cdot))=\ln\alpha_0
$$
and hence
\begin{equation}
\label{stability-eq2}
\frac{1}{\alpha_0}u^*(t,x)\le u_0(t,x;u_0)\le\alpha_0 u^*(t,x)\quad \forall \,\, t>0,\,\, x\in\mathcal{H}.
\end{equation}

If $\alpha_\infty=1$, then for any $\epsilon>0$, there is $\tau>0$ such that for $t\ge\tau$,
$$
\rho(u(t,\cdot;u_0),u^*(t,\cdot))\le \ln(1+\epsilon).
$$
This implies that
$$
\frac{1}{1+\epsilon}u^*(t,x)\le u(t,x;u_0)\le (1+\epsilon) u^*(t,x)\quad \forall\,\, t\ge \tau,\,\, x\in\mathcal{H}.
$$
Hence
$$
|u(t,x;u_0)-u^*(t,x)|\le \epsilon u^*(t,x)\quad \forall\,\, t\ge \tau,\,\, x\in\mathcal{H}.
$$
It then follows that
$$
\lim_{t\to\infty}\|u(t,\cdot;u_0)-u^*(t,\cdot)\|=0.
$$

Assume $\alpha_\infty>1$. By \eqref{stability-eq2}, there are $\epsilon>0$, $M>0$,  and $\sigma>0$ such that
$$
\epsilon \le u(t,x;u_0)\le M,\,\, \epsilon\le u^*(t,x)\le M\quad \forall\,\, t\ge 0,\,\, x\in\mathcal{H}
$$
and
$$
\rho(u(t,\cdot;u_0),u^*(t,\cdot))\ge \sigma\quad \forall t\ge 0.
$$
By Proposition \ref{basic-part-metric} again, there is $\delta>0$ such that for any $n\ge 1$,
$$
\rho(u(nT,\cdot;u_0),u^*(nT,\cdot))\le \rho(u((n-1)T,\cdot;u_0),u^*((n-1)T,\cdot))-\delta
$$
and hence
$$
\rho(u(nT,\cdot;u_0),u^*(nT,\cdot))\le \rho(u_0,u^*(0,\cdot))-n\delta\quad \forall\,\, n\ge 1.
$$
Let $n\to\infty$, we have
$$\lim_{n\to\infty}\rho(u(nT,\cdot;u_0),u^*(nT,\cdot))=-\infty.
$$
 This is
a contradiction.

Therefore, we must have $\alpha_\infty=1$ and
$$
\lim_{t\to\infty}\|u(t,\cdot;u_0)-u^*(t,\cdot)\|=0.
$$
\end{proof}

\subsection{Existence}

In this subsection, we show the existence of time periodic strictly positive solutions of \eqref{main-eq}, i.e., show Theorem
\ref{positive-solution-thm}(3). To this end, we first prove some lemmas.

 Throughout this subsection, we assume the conditions in
Theorem
\ref{positive-solution-thm}(3). Then by Proposition \ref{basic-periodic-solution}, \eqref{periodic-eq} has a unique
time and space periodic positive solution $u_0^*(t,x)$. Let $\delta_0>0$ be such that
$$
0<\delta_0<\inf_{(t,x)\in\RR\times\mathcal{H}}u_0^*(t,x).
$$
Let $\psi_0:\RR\to\RR^+$ be a non-increasing smooth function such that
\begin{equation}
\label{psi-eq}
0\le \psi_0(\cdot)\le\delta_0, \,\, \,\,\liminf_{r\to -\infty}\psi_0(r)>0,\,\,\,\, \psi_0(r)=0\quad \forall\,\, r\ge 0.
\end{equation}

\begin{lemma}
\label{existence-lm1} For given $\xi\in S^{N-1}$, let $u_0(x)=\psi_0(x\cdot\xi)$ and $u_{n,\xi}(x)=u_0(x+n\xi)$ $(n\in\NN$).
For any $0<c^{'}<c_0^*(\xi)$, there are  $K\ge 0$ and $n^*\ge 0$ such that
$$
u(t+KT,x;u_{n^*,\xi},f)\ge \delta_0
$$
for $t\in [0,KT]$ and $x\in\mathcal{H}$ with $x\cdot\xi\le c^{'}(t+KT)-n^*$, and
$$
u(t+KT,\cdot;u_{n^*,\xi},f)\ge u_{n^*,\xi}(\cdot)\quad \forall \,\,\, t\in [0,KT].
$$
 \end{lemma}

\begin{proof} Let $\epsilon>0$ be such that
$\delta_0<\inf_{(t,x)\in\RR\times\mathcal{H}}u_0^*(t,x)-\epsilon$. Fix $0<c^{'}<c_0^*(\xi)$.
By Proposition \ref{basic-spreading-speed}, there is $K\in \NN$ such that,
\begin{equation}
\label{spread-aux-eq1}
u(t+KT,x;u_0,f_0(\cdot,\cdot+y,\cdot))\geq u_0^{*}(t,x+y)-\epsilon/2
\end{equation}
 for $t\ge 0$, $x\cdot \xi\le c^{'}(t+KT)$
and $y\in\mathcal{H}$.

Observe that
$$
f(t,x-n\xi,u)- f_0(t,x-n\xi,u)\to 0$$
as $n\to\infty$ uniformly in $(t,x,u)$ for $(t,u)$ in bounded sets and $x$ in sets with
$x\cdot\xi$ being bounded. Then by Proposition \ref{basic-convergence},
 there exists
$n^*\in \NN$ such that
\begin{equation}
\label{ineq-1}
u(t+KT,x;u_0,f(\cdot, \cdot-n\xi, \cdot))\geq u(t+KT,x;u_0,f_0(\cdot,\cdot-n\xi,\cdot)-\epsilon/2, \quad n\geq n^*
\end{equation}
for $t\in [0,KT]$ and $x\in\mathcal{H}$ with  $c^{'}(t+KT)-1\le x\cdot\xi\le c^{'}(t+KT)$.
This together with \eqref{spread-aux-eq1} implies that
\begin{equation}
\label{ineq-2}
u(t+KT,x;u_0,f(\cdot, \cdot-n\xi, \cdot))\geq
u_0^{*}(t+KT, x-n\xi)-\epsilon, \quad n\geq n^*
\end{equation}
for $t\in[0,KT]$ and $x\in\mathcal{H}$ with $x\cdot \xi\in [c^{'}(t+KT)-1,c^{'}(t+KT)]$.

Note that $$u(t+KT,x+n\xi;u_0,f(\cdot, \cdot-n\xi, \cdot))=u(t+KT,x;u_0(\cdot+n\xi),f).$$
We then have
 \begin{equation}
 \label{spread-aux-eq2}
 u(t+KT,x;u_0(\cdot+n^*\xi),f)\geq
u_0^{*}(t+KT, x)-\epsilon
\end{equation}
 for $t\in [0,KT]$ and $x\in\mathcal{H}$ with $x\cdot \xi\in [c^{'}(t+KT)-1-n^*,c^{'}(t+KT)-n^*]$.

By Proposition \ref{basic-comparison} and \eqref{ineq-2},
\begin{align*}
u(t+KT,x;u_0(\cdot+n^{*}\xi),f)&=u(t+KT,x+(n^*+1)\xi;u_0(\cdot-\xi),f(\cdot,\cdot-(n^*+1)\xi,\cdot))\\
&\geq u(t+KT,x+(n^*+1)\xi;u_0,f(\cdot,\cdot-(n^*+1)\xi,\cdot))\\
&\geq u_0^{*}(t+KT, x)-\epsilon
\end{align*}
for $t\in [0,KT]$ and $x\in\mathcal{H}$ with $x\cdot \xi\in [c^{'}(t+KT)-(n^*+2),c^{'}(t+KT)-(n^*+1)]$.
This together with \eqref{spread-aux-eq2}
implies that
 $$u(t+KT,x;u_0(\cdot+n^*\xi),f)\geq u_0^{*}(t+KT, x)-\epsilon$$ for
 $t\in [0,KT]$ and $x\in\mathcal{H}$ with
  $x\cdot \xi\in
[c^{'}(t+KT)-(n^*+2), c^{'}(t+KT)-n^{*}]$.

By induction, we have
$$u(t+KT,x;u_0(\cdot+n^*\xi),f)\geq u_0^{*}(t+KT, x)-\epsilon$$ for
$t\in [0,KT]$ and $x\in\mathcal{H}$ with  $x\cdot \xi\in (-\infty,
c^{'}(t+KT)-n^{*}]$. The lemma then follows from the fact that
$u_0(x+n^*\xi)\le \delta_0<\inf_{(t,x)\in\RR\times\mathcal{H}}u_0^*(t,x)-\epsilon$ for all $x\in\mathcal{H}$.
\end{proof}

\begin{lemma}
\label{existence-lm2}
Let $u_0(\cdot)$, $K$, and $n^*$ be as in Lemma \ref{existence-lm1}.
For any $M\gg 1$, there exists $\tilde \delta_0>0$ such that $u(t, x; u_{n^*,\xi},f)\geq \tilde \delta_0$ for $t\ge KT$ and $x\cdot \xi \leq M$.
\end{lemma}

\begin{proof}
By Proposition \ref{basic-comparison} and  Lemma \ref{existence-lm1},
\begin{equation}
\label{spread-auc-eq3}
u(mKT,\cdot;u_{n^*,\xi},f)\ge u_{n^*,\xi}(\cdot)\quad \forall \,\, m\ge 1
\end{equation}
and
\begin{equation}
\label{spread-aux-eq3-1}
u(t,x;u_{n^*,\xi},f)\ge \delta_0\quad \forall\,\, t\ge KT,\,\, x\in\mathcal{H}\,\, {\rm with}\,\,  x\cdot\xi\le -n^*.
\end{equation}
It then suffices to prove that for any $M\gg 1$,
\begin{equation}
\label{spread-aux-eq4}
\inf_{t\ge KT, x\in\mathcal{H}, x\cdot\xi\in [-n^*,M]}u(t,x;u_{n^*,\xi},f)>0.
\end{equation}

Suppose that \eqref{spread-aux-eq4} does not hold. Then there are $t_n\ge KT$ and $x_n\in\mathcal{H}$ with $x_n\cdot\xi\in [-n^*,M]$  such that
\begin{equation}
\label{spread-aux-eq4-1}
u(t_n,x_n;u_{n^*,\xi},f)\to 0\quad {\rm as}\quad n\to\infty.
\end{equation}
Note that there are $k_n\in \NN$ and $\tau_n\in [0,KT]$ such that
$$
t_n=k_nKT+\tau_n.
$$
Then by \eqref{spread-auc-eq3},
$$
u(t_n,x_n;u_{n^*,\xi},f)=u(\tau_n,x_n;u(k_nKT,\cdot;u_{n^*,\xi},f),f)\ge u(\tau_n,x_n;u_{n^*,\xi},f).
$$
This together with \eqref{spread-aux-eq4-1} implies that
\begin{equation}
\label{spread-aux-eq4-2}
u(\tau_n,x_n;u_{n^*,\xi},f)\to 0\quad {\rm as}\quad n\to\infty.
\end{equation}

Without loss of generality, we may assume that $\tau_n\to \tau^*$ as $n\to\infty$ for some
$\tau^*\in [0,KT]$. We then have
$$
\|u(\tau_n,\cdot;u_{n^*,\xi},f)-u(\tau^*,\cdot;u_{n^*,\xi},f)\|\to 0\quad {\rm as}\quad n\to\infty
$$
and hence
$$
u(\tau_n,x_n;u_{n^*,\xi},f)-u(\tau^*,x_n;u_{n^*,\xi},f)\to 0\quad {\rm as}\quad n\to\infty.
$$

In the case that $\{|x_n\|\}$ is unbounded, we may assume  that $\|x_n\|\to \infty$ as $n\to\infty$. Then
$$
f(t,x+x_n,u)-f_0(t,x+x_n,u)\to 0
$$
as $n\to\infty$ uniformly in $(t,x,u)$ on bounded sets.
Observe that
$$
u(\tau^*,x_n;u_{n^*,\xi},f)=u(\tau ^*,0;u_{n^*,\xi}(\cdot+x_n),f(\cdot,\cdot+x_n,\cdot))\quad \forall\,\, n\ge 1.
$$
 Observe also that there is $n^{**}\ge n^*$ such that
$$
u_{n^*,\xi}(\cdot+x_n)\ge u_{n^{**},\xi}(\cdot)\quad \forall\,\, n\ge 1.
$$
By Propositions \ref{basic-comparison} and \ref{basic-convergence}, we have
$$
u(\tau^*,0;u_{n^*,\xi} (\cdot+x_n),f(\cdot,\cdot+x_n,\cdot))\ge u(\tau^*,0;u_{n^{**},\xi}(\cdot),f(\cdot,\cdot+x_n,\cdot))
$$
and
\begin{equation}
\label{spread-aux-eq5}
u(\tau^*,0;u_{n^{**},\xi}(\cdot),f(\cdot,\cdot+x_n,\cdot))-u(\tau^*,0;u_{n^{**},\xi}(\cdot),f_0(\cdot,\cdot+x_n,\cdot))\to 0
\end{equation}
as $n\to\infty$. By the periodicity of $f_0$ in $x$, we may also assume that there is $\tilde x\in\mathcal{H}$ such that
$$
f_0(t,x+x_n,u)\to f_0(t,x+\tilde x,u)\quad {\rm as}\quad n\to\infty
$$
uniformly in $(t,x,u)$ on bounded sets. Then by Proposition \ref{basic-convergence} again,
\begin{equation}
\label{spread-aux-eq6}
u(\tau^*,0;u_{n^{**},\xi}(\cdot),f_0(\cdot,\cdot+x_n,\cdot))\to u(\tau^*,0;u_{n^{**},\xi}(\cdot),f_0(\cdot,\cdot+\tilde x,\cdot))\quad {\rm as}\quad n\to\infty.
\end{equation}
By Proposition \ref{basic-comparison},
\begin{equation}
\label{spread-aux-eq7}
u(\tau^*,0;u_{n^{**},\xi}(\cdot),f_0(\cdot,\cdot+\tilde x,\cdot))>0.
\end{equation}
It then follows from \eqref{spread-aux-eq5}, \eqref{spread-aux-eq6}, and \eqref{spread-aux-eq7} that
$$
\liminf_{n\to\infty}u(\tau_n,x_n;u_{n^*,\xi},f)>0,
$$
which contradicts to \eqref{spread-aux-eq4-2}. Therefore, \eqref{spread-aux-eq4} holds.

In the case that $\{\|x_n\|\}$ is bounded, we may assume that $x_n\to \tilde x\in \mathcal{H}$ as $n\to\infty$.
Then
$$
u(\tau^*,x_n; u_{n^*,\xi},f)\to u(\tau^*,\tilde x;u_{n^*,\xi},f)>0\quad {\rm as}\quad n\to\infty.
$$
This also contradicts to \eqref{spread-aux-eq4-2}. Therefore, \eqref{spread-aux-eq4} holds.
The lemma is thus proved.
\end{proof}

Observe that for any $M\ge M_0$, $u(t,x)\equiv M$ is a supersolution of \eqref{main-eq} on $\mathcal{H}$. Hence
$u(T,x;M,f)\le M$ and then by Proposition \ref{basic-comparison}, $u(nT,x;M,f)$ decreases at $n$ increases.
Define
\begin{equation}
\label{u-plus-eq}
u^{+}( x):={\rm lim}_{n\to \infty}u(nT,x;M,f).
\end{equation}
 Then $u^+(x)$ is a Lebesgue measurable and upper semi-continuous function. In the following, we fix an $M\ge \max\{M_0,\delta_0\}$.

\begin{lemma}
\label{lm-lift2}  There exists $\bar \delta >0$ such that $u^+(x)\geq\bar \delta$ for $x\in \mathcal{H}$.
\end{lemma}

\begin{proof}
Let $\psi_0(\cdot)$ be as in \eqref{psi-eq} and
$$
u_{\pm i}(x)=\psi_0(\pm x\cdot {\bf e_i}),\quad i=1,2,\cdots, N.
$$
By Proposition \ref{basic-comparison},
$$
u(t,\cdot;M,f)\ge u(t,\cdot; u_{\pm i},f)
$$
for $t\ge 0$ and $i=1,2,\cdots, N$. By Lemma \ref{existence-lm2}, there are $\bar \delta$ and $\bar T>0$ such that
$$
u(t,x;u_{\pm i},f)\ge \bar \delta\quad \forall\,\, t\ge \bar T,\,\, x\cdot{\bf e_i}\le 0,\,\, i=1,2,\cdots, N.
$$
It then follows that
$$
u(t,x;M,f)\ge\bar\delta\quad \forall\,\,  t\ge \bar T,\,\, x\in\mathcal{H}.
$$
This implies that
$$
u(mT,x;M,f)\geq\bar\delta\quad \forall\,\, m\gg 1,\,\, x\in\mathcal{H}
$$
and then
$$
u^+(x)\ge\bar\delta\quad \forall\,\, x\in\mathcal{H}.
$$
The lemma thus follows.
\end{proof}

Now we prove the existence of time periodic positive solutions

\begin{proof} [Proof of Theorem \ref{positive-solution-thm}(3)]
We first claim that
\begin{equation}
\label{existence-proof-eq1}
\rho(u(nT,\cdot;\bar\delta/2),u(nT,\cdot;M))\to 0
\end{equation}
as $n\to\infty$.
Assume this is not true. Let $\rho_n=\rho(u(nT,\cdot;\bar\delta/2),u(nT,\cdot;M))$ and
$\rho_\infty=\lim_{n\to\infty}\rho_n$ (the existence of this limit follows from
Proposition \ref{basic-part-metric}).
Then $\rho_\infty>0$,
$$
\frac{\bar\delta}{e^{\rho_0}}\le\frac{1}{e^{\rho_0}}u(nT,\cdot;M)\le u(nT,\cdot;\bar\delta/2)\le e^{\rho_0}u(nT,\cdot;M)\le e^{\rho_0}M
$$
for $n=0,1,2,\cdots$ and
$$
\rho(u(nT,\cdot;\bar\delta/2),u(nT,\cdot;M))\ge \rho_\infty
$$
for $n=1,2,\cdots$. By Proposition \ref{basic-part-metric}, there is $\delta>0$ such that
$$
\rho(u(nT,\cdot;\bar\delta/2),u(nT,\cdot;M))\le \rho_0-n\delta\quad \forall\,\, n=1,2,\cdots.
$$
This implies that $\rho_\infty=-\infty$, a contradiction. Therefore, \eqref{existence-proof-eq1} holds.

By \eqref{existence-proof-eq1}, there is $K_1\ge 1$ such that
$$
\bar\delta/2\le u(K_1T,\cdot;\bar\delta/2)
$$
and then
$$
\bar\delta/2\le u(nK_1T,\cdot;\bar\delta/2)\le u(nK_1T,\cdot;M)\le M\quad \forall\,\, n=1,2,\cdots.
$$
It then follows that
$$
u(nK_1T,x;M)\ge u^+(x)\ge u(nK_1T,x;\bar\delta/2)\quad \forall\,\, x\in\mathcal{H},\,\, n=1,2,\cdots.
$$
Therefore
\begin{align*}
0\le u(nK_1T,x;M)-u^+(x)&\le u(nK_1T,x;M)-u(nK_1T,x;\bar\delta/2)\\
&\le u(nK_1T,x;M)(1-\frac{1}{e^{\rho_n}})\le M(1-\frac{1}{e^{\rho_n}}).
\end{align*}
This implies that
$$
\lim_{n\to\infty}u(nK_1T,x;M)=u^+(x)
$$
uniformly in $x\in\mathcal{H}$ and $u^+(\cdot)\in X^{++}$. Moreover,
by
$$
u(nK_1T,\cdot;M)\ge u(kT,\cdot;M)\ge u((n+1)K_1,\cdot;M)\quad \forall \,\, nK_1\le k\le (n+1)K_1,
$$
we have
$$
\lim_{k\to\infty}u(kT,x;M)=u^+(x)
$$
uniformly in $x\in\mathcal{H}$ and then
$$
u(T,\cdot;u^+)=u^+(\cdot).
$$
This implies that $u^*(t,x)=u(t,x;u^+)$ is a time periodic strictly positive solutions of \eqref{main-eq}.
\end{proof}

\subsection{Tail property}

In this subsection, we prove the tail property of time periodic strictly positive solutions of \eqref{main-eq}.
Throughout this subsection, we assume the conditions in  Theorem \ref{positive-solution-thm} (4).

\begin{proof} [Proof of Theorem \ref{positive-solution-thm} (4)]
Suppose that $u^*(t,x)$ is a time periodic strictly positive solution of \eqref{main-eq}.
Observe that $u^*(t,x)=u(t,x;u^+)$, where $u^+$ is as in the proof of Theorem \ref{positive-solution-thm}(3).
We claim
\begin{equation}
\label{tail}
\lim_{r\to\infty}\sup_{x\in\mathcal{H},\|x\|\geq r}|u^*(t, x)-u_0^*(t, x)|=0\quad \forall t\in\RR.
\end{equation}

To prove \eqref{tail}, we first
 show that
\begin{equation}
\label{tails}
\lim_{r\to\infty}\sup_{x\in\mathcal{H},\|x\|\geq r}|u^+( x)-u_0^{+}( x)|=0.
\end{equation}
Recall  $u^{+}( x):={\rm lim}_{n\to \infty}u(nT,x; M,f)$ and  $u_0^{+}( x):={\rm lim}_{n\to \infty}u(nT,x; M, f_0)$.

Assume \eqref{tails} is not true. Then there exists $\epsilon_0>0$ and $\{x_k\}\in \RR$ with $\|x_k\| \to \infty$ such that
$$
|u^+(x_k)-u_0^{+}(x_k)|> \epsilon_0
$$ for $k\ge 1$.

Since both $u(nT,x; M,f)\to u^+(x)$ and $u(nT,x; M,f_0)\to u_0^{+}( x)$ uniformly on $x\in \mathcal{H}$, there is $\tilde N$ such that for $n\geq\tilde  N$,
\begin{equation}
\label{ineq 2}
|u^+(nT, x_k;M,f)-u^{+}(nT, x_k;M, f_0)|> \epsilon_0\quad \forall\,\, k\ge 1.
\end{equation}

Note that there is $\tilde x_0\in\mathcal{H}$ such that
$$
f_0(t,x+x_k,u)\to f_0(t,x+\tilde x_0,u)
$$
as $k\to\infty$ uniformly in $(t,x,u)$ on bounded sets.
Note also that
$$
f(t,x+x_k,u)-f_0(t,x+x_k,u)\to 0
$$
as $k\to\infty$ uniformly in $(t,x,u)$ on bounded sets. Hence
$$
f(t,x+x_k,u)\to f_0(t,x+\tilde x_0,u)
$$
as $k\to\infty$ uniformly in $(t,x,u)$ on bounded sets. Then by Proposition \ref{basic-convergence},
\begin{equation*}
\begin{split}
|u(\tilde NT,x_k;M,f)-u(\tilde NT,x_k;M, f_0)|& = |u(\tilde NT,0;M,f(\cdot,\cdot+x_k; \cdot))-u(\tilde NT,0;M, f_0(\cdot, \cdot+x_k;\cdot))|\\
& \leq |u(\tilde NT,0;M,f(\cdot,\cdot+x_k; \cdot))
-u(\tilde NT,0;M, f_0(\cdot, \cdot+\tilde x_0;\cdot))| \\ &+ |u(\tilde NT,0;M, f_0(\cdot, \cdot+\tilde x_0;\cdot))-u(\tilde NT,0;M, f_0(\cdot, \cdot+x_k;\cdot))|\\
& \to 0
\end{split}
\end{equation*}
as $k\to\infty$.
This contradicts to \eqref{ineq 2}. Therefore, \eqref{tails} holds.

Now we prove \eqref{tail}. Note that we only need to prove \eqref{tail} for $t>0$.  Recall that $u^*(t,x)=u(t,x;u^+,f)$ and $u_0^*(t,x)=u(t,x;u_0^+,f_0)$. Suppose that
\eqref{tail} does not hold for some $t>0$.
Then there are $x_k\in\mathcal{H}$ with $\|x_k\|\to \infty$ and $\epsilon_0>0$ such that
$$
|u(t,x_k;u^+,f)-u(t,x_k;u_0^+,f_0)|\ge \epsilon_0\quad \forall\,\, k\ge 1.
$$
Hence
$$
|u(t,x_k;u^+,f)-u(t,x;u_0^+,f_0)|=|u(t,0;u^+(\cdot+x_k),f(\cdot+x_k))-u(t,0;u_0^+(\cdot+x_k),f_0(\cdot+x_k))|\ge \epsilon_0
$$
for all $k\ge 1$.
By (H1), \eqref{tails}, Proposition \ref{basic-convergence}, and the arguments in the proof of \eqref{tails},
$$
\lim_{k\to\infty} [u(t,0;u^+(\cdot+x_k),f(\cdot+x_k))-u(t,0;u_0^+(\cdot+x_k),f_0(\cdot+x_k))]=0.
$$
This is a contradicts again. Therefore, \eqref{tail} holds.
\end{proof}

\section{Spatial Spreading Speeds}

In this section, we investigate the spatial spreading speeds of
\eqref{main-eq} and prove Theorems \ref{spreading-speed-thm} and \ref{spreading-feature-thm}.
To do so, we first prove a lemma.

Throughout this section, we assume the conditions in Theorem \ref{spreading-speed-thm}.
Let
$u_0^*(t,x)$ be the unique time and space periodic positive solution of \eqref{periodic-eq}. Let $\delta_0>0$ be such that
$$
0<\delta_0<\inf_{(t,x)\in\RR\times\mathcal{H}}u_0^*(t,x).
$$

\begin{lemma}
\label{spread-lm1}
\begin{itemize}
\item[(1)] Let $\xi\in S^{N-1}$, $c>0$ and $u_0\in X^+$ be given.
 If $\liminf_{x\cdot\xi\leq ct,t\to\infty}u(t,x;u_0,f)>0$,
then for any $0<c^{'}<c$,
$$
\limsup_{x\cdot\xi\leq c^{'}t,t\to\infty}|u(t,x;u_0,f)-u^*(t,x)|=0.
$$

\item[(2)]  Let $\xi\in S^{N-1}$, $c>0$ and $u_0\in X^+$ be given.
 If $\liminf_{|x\cdot\xi|\leq ct,t\to\infty}u(t,x;u_0,f)>0$,
then for any $0<c^{'}<c$,
$$
\limsup_{|x\cdot\xi|\leq c^{'}t,t\to\infty}|u(t,x;u_0,f)-u^*(t,x)|=0.
$$

\item[(3)]  Let  $c>0$ and $u_0\in X^+$ be given.
 If $\liminf_{\|x\|\leq ct,t\to\infty}u(t,x;u_0,f)>0$,
then for any $0<c^{'}<c$,
$$
\limsup_{\|x\|\leq c^{'}t,t\to\infty}|u(t,x;u_0,f)-u^*(t,x)|=0.
$$
\end{itemize}
\end{lemma}

\begin{proof}
It can be proved by the similar arguments as in \cite[Lemma 5.1]{KoSh}. For completeness, we provide a proof of (1) in
the following. (2) and (3) can be proved by similar arguments.

(1) Suppose that $\liminf_{x\cdot\xi\leq
ct,t\to\infty}u(t,x;u_0,f)>0$. Then there are $\delta$ and $T^*>0$
such that
\vspace{-0.05in}$$
u(t,x;u_0, f)\geq \delta\quad \forall
(t,x)\in\RR^+\times\mathcal{H},\,\, x\cdot\xi\leq c t,\,\,
t\geq T^*.
\vspace{-0.05in}$$
Assume that the conclusion of (1) is not true. Then there are $0<c^{'}<c$, $\epsilon_0>0$,
$x_n\in\mathcal{H}$, and $t_n\in\RR^+$ with $x_n\cdot\xi\leq c^{'}t_n$
and $t_n\to\infty$ such that
\vspace{-0.05in}\begin{equation}
\label{tech-lm-eq1}
|u(t_n,x_n;u_0,f)-u^*(t_n,x_n)|\geq \epsilon_0\quad \forall n\ge 1.
\vspace{-0.05in}\end{equation}
Note that there are $k_n\in\ZZ^+$  and $\tau_n\in [0,T]$ such that $t_n=k_nT+\tau_n$.
Without loss of generality, we may assume that $\tau_n\to \tau^*$ and  $x_n\to x^*$ as
$n\to\infty$ in the case that $\{\|x_n\|\}$ is bounded (this implies
that $f(t+t_n,x+x_n,u)\to f(t+\tau^*,x+x^*,u)$ uniformly in $(t,x,u)$ on bounded
sets) and $f(t+t_n,x+x_n,u)- f_0(t+\tau^*, x+x_n, u)\to 0$ as $n\to\infty$ uniformly in
$(t,x,u)$ on bounded sets in the case that $\{\|x_n\|\}$ is unbounded.

Let $\tilde u_0\in X^+$, \vspace{-0.05in}$$ \tilde
u_0(x)=\delta\quad \forall x\in\mathcal{H}. \vspace{-0.05in}$$ Let
$$
M=\sup_{x\in\mathcal{H}}u_0(x). $$
 By Theorem
\ref{positive-solution-thm}, there is $\tilde T>0$ such that
 \vspace{-0.05in}\begin{equation} \label{tech-lm-eq2}
|u(t,x;M,f)-u^*(t,x)|<\frac{\epsilon_0}{2}\quad \forall\,\, t\ge
\tilde T,\,\,  x\in\mathcal{H}, \vspace{-0.05in}\end{equation}
\vspace{-0.05in}\begin{equation} \label{tech-lm-eq3} |u(\tilde
T,x;\tilde u_0,f(\cdot+\tau,\cdot+x^*,\cdot))-u^*(\tilde
T+\tau,x+x^*)|<\frac{\epsilon_0}{2},\quad \forall \,
x\in\mathcal{H}, \,\,\tau\in\RR \vspace{-0.05in}\end{equation} and
\vspace{-0.05in}\begin{equation} \label{tech-lm-eq4} |u(\tilde
T,x;\tilde u_0,f_0(\cdot+\tau,\cdot,\cdot))-u_0^*(\tilde
T+\tau,x)|<\frac{\epsilon_0}{2}\quad\forall\, x\in\mathcal{H},\,\,\tau\in\RR.
\vspace{-0.05in}\end{equation}

Without loss of generality, we may assume that $t_n-\tilde T\geq T^*$ for $n\geq 1$.
Let $\tilde u_{0n}\in X^+$ be such that $\tilde u_{0n}(x)=\delta$ for
$x\cdot\xi\leq \frac{c^{'}+c}{2}(t_n-\tilde T)$, $0\le \tilde u_{0n}(x)\le
\delta$ for $\frac{c^{'}+c}{2}(t_n-\tilde T)\le x\cdot\xi\le c
(t_n-\tilde T)$, and $\tilde u_{0n}(x)=0$ for $x\cdot\xi\geq c(t_n-\tilde T)$.
Then
\vspace{-0.05in}$$
u(t_n-\tilde T,\cdot;u_0,f)\geq \tilde u_{0n}(\cdot)
\vspace{-0.05in}$$
and hence
\vspace{-0.05in}\begin{align}
\label{tech-lm-eq5}
u(t_n,x_n;u_0,f)&=u(\tilde T,x_n;u(t_n-\tilde T,\cdot;u_0,f), f(\cdot+ t_n-\tilde T, \cdot,\cdot))\nonumber\\
&=u(\tilde T,0;u(t_n-\tilde T,\cdot+x_n;u_0,f),f(\cdot+t_n-\tilde T,\cdot+x_n,\cdot))\nonumber\\
&\geq u(\tilde T,0;\tilde u_{0n}(\cdot+x_n),f(\cdot+t_n-\tilde T,\cdot+x_n,\cdot)).
\vspace{-0.05in}\end{align}

Observe that $\tilde u_{0n}(x+x_n)\to \tilde u_0(x)$ as $n\to\infty$
uniformly in $x$  on bounded sets. In the case  that
$f(t+t_n,x+x_n,u)- f_0(t+\tau^*,x+x_n,u)\to 0$ as $n\to\infty$,  by Proposition
\ref{basic-convergence}, \vspace{-0.05in}$$ u(\tilde T,0;\tilde
u_{0n}(\cdot+x_n),f(\cdot+t_n-\tilde T, \cdot+x_n,\cdot))- u(\tilde
T,0;\tilde u_0,f_0(\cdot+\tau^*-\tilde T,\cdot+x_n,\cdot))\to 0
\vspace{-0.05in}$$ as $n\to\infty$. Then by \eqref{tech-lm-eq4} and
\eqref{tech-lm-eq5}, \vspace{-0.05in}\begin{equation}
\label{tech-lm-eq6} u(t_n,x_n;u_0,f)>
u_0^*(\tau^*,x_n)-\epsilon_0/2\quad {\rm for}\quad n\gg 1.
\vspace{-0.05in}\end{equation}
By Theorem \ref{positive-solution-thm}(4),
\vspace{-0.05in}\begin{equation} \label{tech-lm-eq7}
u_0^*(\tau^*,x_n)>u^*(\tau^*,x_n)-\epsilon_0/2\quad {\rm for}\quad
n\gg 1. \vspace{-0.05in}\end{equation}
 By Proposition
\ref{basic-comparison} and \eqref{tech-lm-eq2},
$$
u(t_n,x_n;u_0,f)\le u(t_n,x_n;M,f)\le u^*(t_n,x_n)+\epsilon_0\quad
\forall \,\, n\gg 1.
$$
Then by
 \eqref{tech-lm-eq6},
\eqref{tech-lm-eq7}, and the continuity of $u^*(t,x)$,
\vspace{-0.05in}$$ |u(t_n,x_n;u_0,f)-u^*(t_n,x_n)|<\epsilon_0\quad
{\rm for}\quad n\gg 1. \vspace{-0.05in}$$ This contradicts to
\eqref{tech-lm-eq1}.

In the case that  $x_n\to x^*$, by Proposition
\ref{basic-convergence} again, \vspace{-0.05in}$$ u(\tilde
T,0;\tilde u_{0n}(\cdot+x_n),f(\cdot+t_n-\tilde
T,\cdot+x_n,\cdot))\to u(\tilde T,0;\tilde u_0,f(\cdot+\tau^*-\tilde
T, \cdot+x^*,\cdot)) \vspace{-0.05in}$$ as $n\to\infty$. By
\eqref{tech-lm-eq3} and \eqref{tech-lm-eq5},
\vspace{-0.05in}\begin{equation} \label{tech-lm-eq8}
u(t_n,x_n;u_0,f)>u^*(\tau^*,x^*)-\epsilon_0/2\quad {\rm for}\quad
n\gg 1. \vspace{-0.05in}\end{equation} By the continuity of
$u^*(\cdot,\cdot)$, \vspace{-0.05in}\begin{equation}
\label{tech-lm-eq9}
u^*(\tau^*,x^*)>u^*(\tau_n,x_n)-\epsilon_{0}/2\quad {\rm for}\quad
n\gg 1. \vspace{-0.05in}\end{equation} By Proposition
\ref{basic-comparison} and \eqref{tech-lm-eq2},
$$
u(t_n,x_n;u_0,f)\le u(t_n,x_n;M,f)\le u^*(t_n,x_n)+\epsilon_0\quad
\forall \,\, n\gg 1.
$$
This together with \eqref{tech-lm-eq8}, and \eqref{tech-lm-eq9}
implies that
 \vspace{-0.05in}$$
|u(t_n,x_n;u_0,f)-u^*(t_n,x_n)|<\epsilon_0\quad {\rm for}\quad n\gg
1. \vspace{-0.05in}$$ This contradicts to \eqref{tech-lm-eq1} again.

 Hence
\vspace{-0.05in}$$ \limsup_{x\cdot\xi\leq
c^{'}t,t\to\infty}|u(t,x;u_0,f)-u^*(t,x)|=0 \vspace{-0.05in}$$ for
all $0<c^{'}<c$.
\end{proof}

\subsection{Proof of Theorem \ref{spreading-speed-thm}}

In this subsection, we prove Theorem \ref{spreading-speed-thm}.

Observe that for any given $\xi\in S^{N-1}$, there is $\mu^*(\xi)>0$
such that
$$
c_0^*(\xi)=\frac{\lambda_{\xi,\mu^*(\xi)}(f_0(\cdot,\cdot,0))}{\mu^*(\xi)}
$$
and
$$
\frac{\lambda_{\xi,\mu}(f_0(\cdot,\cdot,0))}{\mu}>c_0^*(\xi)\quad
\forall\,\, 0<\mu<\mu^*(\xi).
$$

\begin{proof}[Proof of Theorem \ref{spreading-speed-thm}]
We first show that for any $c>c_0^*(\xi)$,
\begin{equation}
\label{spread-thm-eq1} \limsup_{x\cdot\xi\ge ct,t\to\infty}
u(t,x;u_0,f)=0.
\end{equation}
Let $0<\tilde\mu <\mu<\mu^*(\xi)$ be such that
$$
c=\frac{\lambda_{\xi,\tilde \mu}(f_0(\cdot,\cdot,0))}{\tilde \mu}>\frac{\lambda_{\xi,\mu}(f_0(\cdot,\cdot,0))}{\mu}>c_0^*(\xi).
$$
By Proposition \ref{basic-criterion}, for any $\epsilon>0$, there are $\delta>0$ and
$a(\cdot,\cdot)\in\mathcal{X}_p$ such that
$$
a(t,x)\ge f_0(t,x,0)+\delta,
$$
$\lambda_{\xi,\mu}(a)$ is the principal eigenvalue of $-\p_t+\mathcal{A}_{\xi,\mu}+a(\cdot,\cdot)\mathcal{I}$
with a positive principal eigenfunction $\phi(\cdot,\cdot)\in\mathcal{X}_p$, and
$$
c_0^*(\xi)<\frac{\lambda_{\xi,\mu}(a)}{\mu}<c.
$$
 Let
$$
u_M(t,x)=Me^{-\mu(x\cdot\xi-\frac{\lambda_{\xi,\mu}(a)}{\mu}t)}\phi(t,x).
$$
It is not difficult to verify that $u_M(t,x)$ is a solution of
\begin{equation}
\label{spread-thm-eq2}
u_t=\mathcal{A}u+a(t,x)u,\quad x\in\mathcal{H}.
\end{equation}
Observe that
$$
f_0(t,x,u)\le f_0(t,x,0)\le a(t,x)-\delta\quad \forall\,\, t\in\RR,\,\, x\in\mathcal{H},\,\, u\ge 0.
$$
This together with (H1) implies that there is $\tilde M>0$ such that
\begin{equation}
\label{spread-thm-eq3}
f(t,x,u)\le f(t,x,0)\le a(t,x)\quad \forall \,\, t\in\RR,\,\, x\cdot\xi\ge \tilde M,\,\, u\ge 0.
\end{equation}
Let $M\ge \|u_0\|$ be such that
\begin{equation}
\label{spread-thm-eq4}
f(t,x,u_M(t,x))\le a(t,x)\quad \forall\,\, t\ge 0,\,\, x\cdot\xi\le\tilde M.
\end{equation}
It then follows from \eqref{spread-thm-eq2}, \eqref{spread-thm-eq3}, and \eqref{spread-thm-eq4} that
$u_M(t,x)$ is a super-solution of \eqref{main-eq}. By Proposition \ref{basic-comparison},
$$
u(t,x;u_0,f)\le u_M(t,x)=Me^{-\mu(x\cdot\xi-\frac{\lambda_{\xi,\mu}(a)}{\mu}t)}\phi(t,x)\quad\forall t\ge 0,\,\, x\in\mathcal{H}.
$$
This implies that \eqref{spread-thm-eq1} holds.

Next, we prove that for any $c<c_0^*(\xi)$,
\begin{equation}
\label{spread-thm-eq5}
\limsup_{x\cdot\xi\le ct,t\to\infty}|u(t,x;u_0,f)-u^*(t,x)|=0.
\end{equation}
First of all, by Proposition \ref{basic-criterion}, there is $\epsilon>0$ such that
$$
c<\inf_{\mu>0} \frac{\lambda_{\xi,\mu}(f_0(\cdot,\cdot,0)-\epsilon)}{\mu}<c_0^*(\xi).
$$
By (H1),  there is $M>0$ such that
\begin{equation}
\label{spread-thm-eq6}
f(t,x,u)\ge f_0(t,x,u)-\epsilon\quad \forall \, t\in\RR, \,\, x\cdot\xi \ge M,\,\, 0\le u\le 1.
\end{equation}
By Lemma \ref{existence-lm2}, there are $\tilde\delta>0$ and $\tilde T>0$ such that
\begin{equation}
\label{spread-thm-eq7}
u(t,x;u_0,f)\ge \tilde \delta\quad \forall \,\, t\ge \tilde T,\,\, x\cdot\xi\le M.
\end{equation}

For given $K^*>0$,
 consider equation
 \begin{equation}
 \label{spread-thm-eq8}
 u_t=(\mathcal{A} u)(t,x)+[f_0(t,x,u)-\epsilon-K^*u]u(x),\quad x\in\mathcal{H}.
\end{equation}
By Proposition \ref{basic-periodic-solution}, \eqref{spread-thm-eq8} has a unique time and space periodic
solution $u_{0,K^*}^*(t,x)$. Let $K^*\gg 1$ be such that
$$
u_{0,K^*}(\tilde T,x)\le \tilde\delta.
$$
Let $\tilde u_0\in X^+(\xi)$ be such that
$$
\tilde u_0(x)\le \min\{u_{0,K^*}^*(\tilde T,x),u(\tilde T,x;u_0,f)\}.
$$
By Proposition \ref{basic-comparison},
$$
u(t+\tilde T,x;u_0,f)\ge u(t,x;\tilde u_0,\tilde f_0(\cdot+\tilde T,\cdot,\cdot))\quad \forall\,\, t\ge 0,\,\, x\cdot\xi\ge M,
$$
where  $\tilde f_0(t,x,u)=f_0(t,x,u)-\epsilon-K^*u$. By Proposition \ref{basic-spreading-speed}, for any $c<c^{'}<\inf_{\mu>0}\frac{\lambda_{\xi,\mu}(f_0-\epsilon)}{\mu}$,
\begin{equation}
\label{spread-thm-eq9}
\limsup_{x\cdot\xi\le c^{'}t,t\to\infty}|u(t,x;\tilde u_0,\tilde f_0(\cdot+\tilde  T,\cdot,\cdot))-u_{0,K^*}^*(t+\tilde T,x)|=0.
\end{equation}
By \eqref{spread-thm-eq7} and \eqref{spread-thm-eq9},
$$
\liminf_{x\cdot\xi\le c^{'}t,t\to\infty}u(t,x;u_0,f)>0.
$$
This together with Lemma \ref{spread-lm1} (1) implies that \eqref{spread-thm-eq5} holds.
\end{proof}

\subsection{Proof of Theorem \ref{spreading-feature-thm}}

In this subsection, we give a proof of  Theorem \ref{spreading-feature-thm}.

\begin{proof}[Proof of Theorem \ref{spreading-feature-thm}]

(1) For given $u_0\in X^+$ satisfying the conditions in Theorem  \ref{spreading-feature-thm} (1), there are $u_0^\pm\in X^+(\pm\xi)$ such that
$$
u_0(\cdot)\leq u_0^\pm(\cdot).
$$
Then for every $c>\max\{c^*(\xi),c^*(-\xi)\}$,
$$
\limsup_{x\cdot\xi\leq ct, t\to\infty}u(t,x;u_0^+)=0,\quad \limsup_{x\cdot (-\xi)\leq
ct,t\to\infty}u(t,x;u_0^-)=0.
$$
 By Proposition \ref{basic-comparison},
$$
u(t,x;u_0)\leq u(t,x;u_0^\pm)\quad {\rm for}\quad t\geq 0,\, x\in\RR^N.
$$
It then follows that
$$
\limsup_{|x\cdot\xi|\geq ct,t\to\infty}u(t,x;u_0)=0.
$$

(2) First, we show that for any $\sigma>0$, there is $r_\sigma>0$ such that for any $r\ge r_\sigma$,  $u_0\in X^+$ with
$u_0(x)\ge \sigma$ for $|x\cdot\xi|\le r$,
\begin{equation}
\label{spread-feature-eq}
\limsup_{|x\cdot\xi|\le ct, t\to\infty}|u(t,x;u_0,f)-u^*(t,x)|=0.\quad \forall \,\, c<\inf\{c^*(\xi),c^*(-\xi)\}.
\end{equation}
Note that we only need to consider $\sigma$
satisfying $0<\sigma<\min\{\inf u_0^*,\inf u^*\}$.

Let $\xi\in S^{N-1}$ be given.  For given  $0<\sigma< \min\{\inf u_0^*,\inf u^*\}$, let $\tilde v^\sigma(\cdot)\in
C(\RR^+,\RR)$ be such that $\tilde v^\sigma(r)\geq 0$ for $r\in\RR^+$ and
$$
\tilde v^\sigma(r)=\begin{cases} \sigma ,\quad 0\le r\le 1\cr 0,\quad r\geq 2.
\end{cases}
$$
Let
$$
\tilde u^{\sigma,\xi}(x)=\tilde v^\sigma(|x\cdot\xi|).
$$
Then by Proposition \ref{basic-spreading-feature},
 for any $0<c<\min\{c^*(\xi),c^*(-\xi)\}$,
$$
\limsup_{|x\cdot\xi|\le ct,t\to\infty}|u(t,x;\tilde u^{\sigma,\xi},f_0(\cdot,\cdot+z,\cdot))-u_0^*(t,x+z)|=0
$$
uniformly in $z\in\RR^N$.
Hence for any $0<c^{'}<c$, there are  $K\in\NN$ and $\tilde c\in (c^{'},c)$ such that
\begin{equation}
\label{spread-feature-eq0}
c^{'}KT< \tilde c KT< cKT-2,\quad \tilde cKT\in\NN
\end{equation}
 and
\begin{equation}
\label{spread-feature-eq00}
u(t+KT,x;\tilde u^{\sigma,\xi},f_0(\cdot,\cdot+z,\cdot))>\sigma\quad \forall\,\, t\in [0,KT],\,\, x\in\mathcal{H}\,\, {\rm with}\,\,  |x\cdot\xi|\le c(t+KT),\,\, z\in\RR^N.
\end{equation}

We claim that there is $n^*>2$ such that
\begin{equation}
\label{spread-feature-eq000}
\liminf_{|x\cdot\xi|\le \tilde ct,t\to\infty }u(t,x;\tilde u_0,f)>0,
\end{equation}
where
$$
\tilde u_0(x)=\begin{cases}\tilde u^{\sigma,\xi}(x+n^*\xi)\quad {\rm  for}\quad |x\cdot\xi+ n^*|\le 2\cr\cr
\tilde u^{\sigma,\xi}(x-n^*\xi)\quad {\rm  for}\quad  |x\cdot\xi-n^*|\le 2\cr\cr
0\quad {\rm otherwise}
\end{cases}
$$

To prove the above claim, we first
note that
$$
\lim_{n\to\infty} [f(t,x\pm n\xi,u)- f_0(t,x\pm n\xi,u)]=0
$$
uniformly in $(t,x,u)$ on any set $E\subset \RR\times\mathcal{H}\times\RR$ with  $\{(t,x\cdot\xi,u)|(t,x,u)\in E\}$ being a bounded set
(we call such set  $E$ a {\it strip type bounded set}). By Proposition \ref{basic-convergence},
$$
u(t+KT,x;\tilde u^{\sigma,\xi},f(\cdot,\cdot\pm n\xi,\cdot))- u(t+KT,x;\tilde u^{\sigma,\xi},f_0(\cdot,\cdot\pm n\xi,\cdot))\to 0
$$
as $n\to\infty$ uniformly in $(t,x,u)$  on strip type bounded sets.
This together with \eqref{spread-feature-eq00} implies that there is $n^*>2$ such that for any $t\in [0,KT]$,
\begin{equation}
\label{spread-feature-eq1}
u(t+KT,x;\tilde u^{\sigma,\xi},f(\cdot,\cdot\pm n\xi,\cdot))\ge\sigma\quad \forall\,\,  x\in\mathcal{H}\,\, {\rm with}\,\,  |x\cdot\xi|\le c(t+KT),\,\, n\ge n^*.
\end{equation}
Observe that
$$
u(t+KT,x\pm n\xi;\tilde u^{\sigma,\xi}(\cdot\mp n\xi),f)=u(t+KT,x;\tilde u^{\sigma,\xi}(\cdot),f(\cdot,\cdot\pm n\xi,\cdot)).
$$
This together with \eqref{spread-feature-eq1} implies that for any $t\in [0,KT]$ and $n\ge n^*$,
\begin{equation}
\label{spread-feature-eq2}
u(t+KT,x;\tilde u^{\sigma,\xi}(\cdot\mp n\xi),f)\ge \sigma \,\, \forall\,x\in\mathcal{H}\,\,{\rm with}\,\,
-c(t+KT)\pm n\le x\cdot\xi \le c(t+KT)\pm n.
\end{equation}

Next,
by \eqref{spread-feature-eq2}, we have
\begin{align}
\label{spread-feature-eq3}
u(t+2KT,x;\tilde u^{\sigma,\xi}(\cdot\mp n^*\xi),f)&=u(t+KT,x;u(KT,\cdot;\tilde u^{\sigma,\xi}(\cdot\mp n^*\xi),f),f)\nonumber\\
&\ge u(t+KT,x;\tilde u^{\sigma,\xi}(\cdot\mp (n^*+i)\xi,f)
\end{align}
for $t\in [0,KT]$ and  $0\le i\le \tilde cKT$.
It then follows from \eqref{spread-feature-eq2} and \eqref{spread-feature-eq3} that for any $t\in [0,KT]$,
$$
\begin{cases}
u(t+2KT,x;\tilde u^{\sigma,\xi}(\cdot- n^*\xi),f)\ge \sigma\cr
\qquad\qquad\qquad  \forall\,x\in\mathcal{H}\,\, {\rm with}\,\,  -c(t+KT)+n^*\le x\cdot\xi\le c(t+KT)+\tilde cKT+n^*,\cr\cr
u(t+2KT,x;\tilde u^{\sigma,\xi}(\cdot+ n^*\xi),f)\ge \sigma\cr
\qquad\qquad\qquad \forall \, x\in\mathcal{H}\,\, {\rm with}\,\, -c(t+KT)-\tilde cKT-n^*\le x\cdot\xi\le c(t+KT)-n^*.
\end{cases}
$$
By induction, we have that for any $t\in [0,KT]$ and $n\ge 1$,
\begin{equation}
\label{spread-feature-eq4}
\begin{cases}
u(t+nKT,x;\tilde u^{\sigma,\xi}(\cdot- n^*\xi),f)\ge \sigma\cr
\qquad\qquad  \forall\, x\in\mathcal{H}\,\, {\rm with}\,\,  -c(t+KT)+n^*\le x\cdot\xi\le c(t+KT)+ (n-1)\tilde cKT+n^*,\cr\cr
u(t+nKT,x;\tilde u^{\sigma,\xi}(\cdot+ n^*\xi),f)\ge \sigma\cr
\qquad\qquad \forall\, x\in\mathcal{H}\,\,{\rm with}\,\, -c(t+KT)-(n-1)\tilde cKT-n^*\le x\cdot\xi\le c(t+KT)-n^*.
\end{cases}
\end{equation}

Now we show that  $\inf_{t\in [0,KT], x\in\RR^N, |x\cdot\xi|\le n^*-2} u(t+KT,x;\tilde u_0,f)>0$.
Assume that there are $t_n\in [0,KT]$,  $\tilde x_n\in\mathcal{H}$ with $|\tilde x_n\cdot\xi|\le n^*-2$ such that
\begin{equation}
\label{spread-feature-eq5}
u(t_n+KT,\tilde x_n;\tilde u_0,f)\to 0\quad {\rm as}\quad n\to\infty.
\end{equation}
Observe that for any $n\ge 1$, there are $x_n\in {\rm span}\{\xi\}$,  $\eta_n\in\RR^N$, and $R>0$  such that
$$
\tilde x_n=x_n+\eta_n,\,\,\, \|x_n\|\le R,\,\,  \eta_n\cdot\xi=0.
$$
Observe also that
\begin{equation*}
u(t_n+KT,\tilde x_n;\tilde u_0,f)=u(t_n+KT,x_n;\tilde u_0(\cdot+\eta_n),f(\cdot,\cdot+\eta_n,\cdot))=u(t_n+KT,x_n;\tilde u_0,f(\cdot,\cdot+\eta_n,\cdot))
\end{equation*}
and
$$
\lim_{n\to\infty} [f(t,x+\eta_n,u)-f_0(t,x+\eta_n,u)]=0
$$
uniformly in $(t,x,u)$ on bounded set. Without loss of generality, we may assume that there are $t^*\in[0,KT]$ and  $\eta^*\in\mathcal{H}$ such that
$$
\lim_{n\to\infty} t_n=t^*,\,\,\, \lim_{n\to\infty}f_0(t,x+\eta_n,u)=f_0(t,x+\eta^*,u)
$$
uniformly in $(t,x,u)$ on bounded set.  By Proposition \ref{basic-comparison},
$$
\inf_{x\in\mathcal{H},\|x\|\le R}u(t^*+KT,x;\tilde u_0,f_0(\cdot,\cdot+\eta^*,\cdot))>0.
$$
This together with Proposition \ref{basic-convergence} implies that
$$
\liminf_{n\to\infty} u(t_n+KT,\tilde x_n;\tilde u_0,f)>0.
$$
This is a contradiction.
Therefore, we have
\begin{equation}
\label{spread-feature-eq5-1}
\inf_{t\in [0,KT], x\in\mathcal{H}, |x\cdot\xi|\le n^*-2} u(t+KT,x;\tilde u_0,f)>0.
\end{equation}

By \eqref{spread-feature-eq4} and \eqref{spread-feature-eq5-1},
$$
u(KT,\cdot;\tilde u_0,f)\ge \tilde u_0(\cdot).
$$
This together with Proposition \ref{basic-comparison} implies that for $n\ge 2$,
\begin{equation}
\label{spread-feature-eq5-2}
u(t+nKT,x;\tilde u_0,f)=u(t+KT,x;u((n-1)KT,\cdot;\tilde u_0,f),f)\ge u(t+KT,x;\tilde u_0,f).
\end{equation}
\eqref{spread-feature-eq000} then follows from \eqref{spread-feature-eq4}, \eqref{spread-feature-eq5-1}, and \eqref{spread-feature-eq5-2}.

For given $0<\sigma<\min\{\inf u_0^*,\inf u^*\}$, let
$$
r_\sigma=n^*+2.
$$
For any $r\ge r_\sigma$ and  $u_0\in X^+$ with $u_0(x)\ge \sigma$ for $|x\cdot\xi|\le r$,
we have
$$
u_0(\cdot)\ge \tilde u_0(\cdot).
$$
By \eqref{spread-feature-eq000} and Proposition \ref{basic-comparison},
$$
\liminf_{|x\cdot\xi|\le \tilde ct, t\to\infty} u(t,x;u_0,f)>0.
$$
Then by Lemma \ref{spread-lm1} (2),
$$
\limsup_{|x\cdot\xi|\le c^{'}t,t\to\infty}|u(t,x;u_0,f)-u^*(t,x)|=0
$$
for any $0<c^{'}<c<\min\{c^*(\xi),c^*(-\xi)\}$, and hence
$$
\limsup_{|x\cdot\xi|\le ct,t\to\infty}|u(t,x;u_0,f)-u^*(t,x)|=0
$$
for any $0<c<\min\{c^*(\xi),c^*(-\xi)\}$.

Now for any $0<c<\min\{c^*(\xi),c^*(-\xi)\}$,   $\sigma>0$,  and $r>0$, suppose that
 $u_0\in X^+$ satisfies $u_0(x)\geq \sigma$ for all $x\in\mathcal{H}$ with
$|x\cdot \xi|\leq r$. Note that there is $m>0$ such that
$$
-1+f(t,x,u(t,x;u_0))\geq -m \quad \forall x\in\RR^N,\,\, t\geq 0.
$$
Then
$$
u_t(t,x;u_0)\geq \int_{\RR^N}k(y-x)u(t,y;u_0)dy-m u(t,x;u_0)
$$
and hence
$$
(e^{mt}u(t,x;u_0))_t\geq \int_{\RR^N} k(y-x)e^{mt}u(t,y;u_0)dy.
$$
This together with Proposition \ref{basic-comparison} implies that
$$
e^{mt}u(t,\cdot;u_0)\geq e^{\mathcal{K}t}u_0,
$$
where $e^{\mathcal{K}t}=I+\mathcal{K}t+\frac{\mathcal{K}^2 t^2}{2!}+\cdots$ and $\mathcal{K}u$ is defined as in \eqref{k-op}. It is then not difficult to see that there is
$\rho\in (0,1)$ such that
$$\rho\sigma<\inf_{x\in\RR^N}u^*(t,x)\quad {\rm and}\quad
u(T,x;u_0)\geq \rho \sigma\quad {\rm for}\quad |x\cdot\xi|\leq r_\sigma.
$$
Let $v_0(x)=\frac{1}{\rho} u(T,x;u_0)$. By \eqref{spread-feature-eq},
\begin{equation}
\label{e-eq4}
 \limsup_{|x\cdot \xi|\leq ct, t\to\infty}|u(t,x;v_0)-u^*(t,x)|=0.
\end{equation}
By Proposition \ref{basic-comparison} and (H0), we have
\begin{equation}
\label{e-eq5}
u(T+t,x;u_0)\equiv u(t,x;\rho v_0)\geq \rho u(t,x;v_0).
\end{equation}
By \eqref{e-eq4} and \eqref{e-eq5}, there is $\tilde n>0$ such that
\begin{equation}
\label{e-eq6}
u(\tilde nT,x;u_0)\geq \rho\sigma \quad {\rm for}\quad |x\cdot\xi|\leq r_{\rho\sigma}.
\end{equation}
By \eqref{spread-feature-eq} again,
\begin{equation}
 \label{e-eq7}
 \limsup_{|x\cdot\xi|\le ct, t\to\infty}|u(t+T,x;u_0)-u^*(t,x)|=0.
\end{equation}
(2) then follows.

(3)  Fix $c>\sup_{\xi\in S^{N-1}}c^*(\xi)$.

 First, let $u_0$ be as in Theorem \ref{spreading-feature-thm} (3).
  For every given $\xi\in S^{N-1}$, there is $\tilde u_0(\cdot;\xi)\in X^+(\xi)$ such that
$u_0(\cdot)\leq \tilde u_0(\cdot;\xi)$. Then by Proposition \ref{basic-comparison},
$$
0\leq u(t,x;u_0)\leq u(t,x;\tilde u_0(\cdot;\xi))
$$
for $t>0$ and $x\in\mathcal{H}$. It then follows that for any $\xi\in S^{N-1}$,
$$
0\leq \limsup_{x\cdot\xi\geq ct,t\to\infty}u(t,x;u_0)\leq \limsup_{x\cdot\xi\geq ct,t\to\infty}u(t,x;\tilde
u_0(\cdot;\xi))=0.
$$

Take any $c^{'}>c$. Consider all $x\in\RR^N$ with $\|x\|=c^{'}$. By the compactness of
$\p B(0,c^{'})=\{x\in\RR^N|\,\|x\|=c^{'}\}$, there are $\xi_1,\xi_2,\cdots,\xi_K\in S^{N-1}$ such that for every
$x\in \p B(0,c^{'})$, there is $k$ ($1\leq k\leq K$) such that $x\cdot \xi_k\geq c$. Hence for every $x\in\RR^N$ with
$\|x\|\geq c^{'}t$, there is $1\leq k\leq K$ such that $x\cdot
\xi_k=\frac{\|x\|}{c^{'}}\Bigl(\frac{c^{'}}{\|x\|}x\Bigr)\cdot \xi_k\geq \frac{\|x\|}{c^{'}}c\geq ct$. By the
above arguments,
$$
0\leq \limsup_{x\cdot\xi_k\geq ct,t\to\infty}u(t,x;u_0)\leq \limsup_{x\cdot\xi_k\geq
ct,t\to\infty}u(t,x;\tilde u_0(\cdot;\xi_k))=0
$$
for $k=1,2,\cdots K$.
 This implies that
$$
\limsup_{\|x\|\geq c^{'}t, t\to\infty}u(t,x;u_0)=0.
$$
 Since $c^{'}>c$ and $c>\sup_{\xi\in
S^{N-1}}c^*(\xi)$ are arbitrary, we have that for $c>\sup_{\xi\in S^{N-1}}c_{\sup}^*(\xi)$,
$$
\limsup_{\|x\|\geq ct, t\to\infty}u(t,x;u_0)=0.
$$

(4)
As in (2), we first prove that for given
  $0<\sigma<\min\{\inf u_0^*,\inf u^*\}$, there is $r_\sigma>0$ such that  for any $0<c<\inf\{c^*(\xi)|\xi\in S^{N-1}\}$,
  $r\ge r_\sigma$, and  $u_0\in X^+$ with $u_0(x)\ge \sigma$ for
  $\|x\|\le r$,
  \begin{equation}
  \label{new-spread-eq}
  \limsup_{\|x\|\le ct,t\to\infty}|u(t,x;u_0,f)-u^*(t,x)|=0.
  \end{equation}

  To this end, for given
  $0<\sigma<\min\{\inf u_0^*,\inf u^*\}$,
  let
$v_0^\sigma\in C(\RR^+,\RR)$ be such that $v_0^\sigma(r)\ge 0$ for $r\in\RR^+$,
  $v_0^\sigma(r)=\sigma$ for $0\le r\le 1$, and
$v_0^\sigma(r)=0$ for $r\ge 2$. Let $\tilde u_0^\sigma(\cdot)\in X^+$ be such that
$$
\tilde u_0^\sigma(x)=v_0^\sigma(\|x\|).
$$
By Proposition \ref{basic-spreading-feature},  for any $0<c<\inf\{c^*(\xi)|\xi\in S^{N-1}\}$,
$$
\lim_{\|x\|\le ct, t\to\infty}|u(t,x;\tilde u_0^\sigma,f_0(\cdot,\cdot+z,\cdot))-u_0^*(t,x+z)|=0
$$
uniformly in $z\in\RR^N$. This implies that for any $0<c^{'}<c$, there are  $K\in\NN$ and $\tilde c\in (c^{'},c)$ such that
\begin{equation}
\label{new-spread-eq1}
c^{'}KT<\tilde cKT<cKT-3,\quad \tilde c KT\in \NN,
\end{equation}
and
 for any $t\in [0,KT]$,
\begin{equation}
\label{new-spread-eq2}
u(t+KT,x;\tilde u_0^\sigma,f_0(\cdot,\cdot+z,\cdot))>\sigma\quad \forall\,\, x\in\mathcal{H}\,\, {\rm with}\,\, \|x\|\le c(t+KT),\,\, z\in\mathcal{H}.
\end{equation}

We clam that there is $n^*>3$ such that
\begin{equation}
\label{new-spread-eq3}
\liminf_{\|x\|\le \tilde ct,t\to\infty} u(t,x;\tilde u_0,f)>0,
\end{equation}
where
$$
\tilde u_0(x)=\begin{cases}
v_0^\sigma(\|x-n^*\xi\|)\quad {\rm if}\quad x\in{\rm span}\{\xi\} \,\, \text{for some}\,\, \xi\in S^{N-1},\,\,   n^*-2\le \|x\|\le n^*+2\cr\cr
0\quad {\rm for}\quad \|x\|<n^*-2\quad {\rm or}\quad \|x\|>n^*+2
 \end{cases}
 $$

To prove the claim, we first note that
$$
u(t,x;\tilde u_0^\sigma(\cdot-n\xi),f)=u(t,x-n\xi;\tilde u_0^\sigma(\cdot),f(\cdot,\cdot+n\xi,\cdot))
$$
and
$$
\lim_{n\to\infty} [f(t,x-n\xi,u)-f_0(t,x-n\xi,u)]=0
$$
uniformly in $(t,x,\xi,u)$ on bounded sets of $\RR\times\RR^N\times  S^{N-1}\times \RR$. By \eqref{new-spread-eq2},
there is $n^*>3$ such that for any $t\in [0,KT]$ and $\xi\in S^{N-1}$,
\begin{equation}
\label{new-spread-eq4}
u(t+KT,x;\tilde u_0^\sigma(\cdot-n\xi),f)>\sigma\quad \forall \,\, x\in\mathcal{H}\,\, {\rm with}\,\, \|x-n\xi\|\le c(t+KT),\,\, n\ge n^*.
\end{equation}

Next, note that by \eqref{new-spread-eq4},
$$
u(KT,\cdot;\tilde u_0^\sigma(\cdot-n^*\xi),f)\ge \min\{\tilde u_0^\sigma(\cdot-(n^*+i)\xi)|i=0,1,\cdots, \tilde cKT\}.
$$
This together with \eqref{new-spread-eq4} implies that for any $t\in [0,KT]$, $0\le i\le \tilde cKT$, $\xi\in S^{N-1}$,
$$
u(t+2KT,x;\tilde u_0^\sigma(\cdot-n^*\xi),f)>\sigma \quad \forall\,\, x\in\mathcal{H}\,\, {\rm with}\,\, \|x-(n^*+i)\xi|\le c(t+KT).
$$
By induction, we have that for any $t\in[0,KT]$, $n\ge 1$, $0\le i\le  \tilde c(n-1)KT$, and $\xi\in S^{N-1}$,
\begin{equation}
\label{new-spread-eq5}
u(t+nKT,x;\tilde u_0^\sigma(\cdot-n^*\xi),f)>\sigma\quad\forall\, x\in\mathcal{H}\,\, {\rm with}\,\, \|x-(n^*+i)\xi\|\le c(t+KT).
\end{equation}

Now note that for any $\xi\in S^{N-1}$,
$$
\tilde u_0(\cdot)\ge \tilde u_0^\sigma(\cdot-n^*\xi).
$$
Then by \eqref{new-spread-eq5}, we have
\begin{equation}
\label{new-spread-eq6}
u(t+nKT,x;\tilde u_0,f)>0\quad \forall \,\, x\in\mathcal{H}\,\, {\rm with}\,\, c(t+KT)-n^*\le\|x\|\le c(t+KT)+\tilde c(n-1)KT+n^*.
\end{equation}
It is clear that
$$
\inf_{\|x\|\le 2cKT-n^*,t\in [0,KT]}u(t+KT,x;\tilde u_0,f)>0
$$
and then
$$
u(KT,\cdot;\tilde u_0,f)>\tilde u_0.
$$
By Proposition \ref{basic-comparison}, we have
\begin{equation}
\label{new-spread-eq7}
\inf_{t\ge KT,\|x\|\le 2cKT-n^*}u(t,x;\tilde u_0,f)>0.
\end{equation}
\eqref{new-spread-eq3} then follows from \eqref{new-spread-eq6} and \eqref{new-spread-eq7}.

Let now
$$
r_\sigma=n^*+3.
$$
Then for any $r\ge r_\sigma$ and $u_0\in X^+$ with $u_0(x)\ge \sigma$ for $\|x\|\le r$,
$$
u_0(\cdot)\ge \tilde u_0(\cdot).
$$
By \eqref{new-spread-eq3},
$$
\liminf_{\|x\|\le\tilde ct,t\to\infty} u(t,x;u_0,f)>0.
$$
By Lemma \ref{spread-lm1},
$$
\limsup_{\|x\|\le c^{'}t,t\to\infty}|u(t,x;u_0,f)-u^*(t,x)|=0
$$
for any $0<c^{'}<c<\inf\{c^*(\xi)|\xi\in S^{N-1}\}$, which implies \eqref{new-spread-eq} holds.

Finally, by the similar arguments as those in (2), for any $\sigma>0$, $r>0$, and $u_0\in X^+$ with
$u_0(x)\ge \sigma$ for $\|x\|\le r$,
$$
\liminf_{\|x\|\le ct,t\to\infty} |u(t,x;u_0,f)-u^*(t,x)|=0
$$
for any $0<c<\inf\{c^*(\xi)|\xi\in S^{N-1}\}$. (4) is thus proved.
\end{proof}

\end{document}